\documentclass[british,11pt]{extarticle}

\usepackage[T1]{fontenc}
\usepackage[latin9]{inputenc}
\usepackage{color}
\usepackage{babel}
\usepackage{url}
\usepackage{amsmath}
\usepackage{amsthm}
\usepackage{amssymb}
\usepackage{graphicx}
\usepackage[a4paper]{geometry}
\usepackage[pdfusetitle,
 bookmarks=true,bookmarksnumbered=false,bookmarksopen=false,
 breaklinks=false,pdfborder={0 0 1},backref=false,colorlinks=true]
 {hyperref}

\makeatletter
\newcommand{\lyxaddress}[1]{
	\par {\raggedright #1
	\vspace{1.4em}
	\noindent\par}
}
\theoremstyle{plain}
\newtheorem{thm}{\protect\theoremname}
\theoremstyle{remark}
\newtheorem{rem}[thm]{\protect\remarkname}
\theoremstyle{definition}
\newtheorem{defn}[thm]{\protect\definitionname}

\usepackage{tikz}
\usepackage{tikz-cd}
\usetikzlibrary{decorations.pathmorphing}
\usepackage{circledsteps}
\newcommand\leftidx[2]{{\vphantom{#2}}#1#2}

\makeatother

\usepackage[style=numeric]{biblatex}
\providecommand{\definitionname}{Definition}
\providecommand{\remarkname}{Remark}
\providecommand{\theoremname}{Theorem}

\begin{document}
\title{Data-driven modelling of autonomous and forced dynamical systems}
\author{Robert Szalai}
\maketitle

\lyxaddress{School of Engineering Mathematics and Technology, University of Bristol,
Ada Lovelace Building, Tankard's Close, Bristol BS8 1TW, email: r.szalai@bristol.ac.uk}
\begin{abstract}
The paper demonstrates that invariant foliations are accurate, data-efficient
and practical tools for data-driven modelling of physical systems.
Invariant foliations can be fitted to data that either fill the phase
space or cluster about an invariant manifold. Invariant foliations
can be fitted to a single trajectory or multiple trajectories. Over
and underfitting are eliminated by appropriately choosing a function
representation and its hyperparameters, such as polynomial orders.
The paper extends invariant foliations to forced and parameter dependent
systems. It is assumed that forcing is provided by a volume preserving
map, and therefore the forcing can be periodic, quasi-periodic or
even chaotic. The method utilises full trajectories, hence it is able
to predict long-term dynamics accurately. We take into account if
a forced system is reducible to an autonomous system about a steady
state, similar to how Floquet theory guarantees reducibility for periodically
forced systems. In order to find an invariant manifold, multiple invariant
foliations are calculated in the neighbourhood of the invariant manifold.
Some of the invariant foliations can be linear, while others nonlinear
but only defined in a small neighbourhood of an invariant manifold,
which reduces the number of parameters to be identified. An invariant
manifold is recovered as the zero level set of one or more of the
foliations. To interpret the results, the identified mathematical
models are transformed to a canonical form and instantaneous frequency
and damping information are calculated.
\end{abstract}

\section{Introduction}

One goal of data-driven modelling is to create mathematical models
that are otherwise impossible to construct from physical laws. The
requirements of such models are that they generalise to unseen data,
are unique to a given physical system and interpretable, so that they
can be used for all things models are use for \cite{epstein2008}. 

\begin{figure}
\begin{centering}
\includegraphics[width=0.6\linewidth]{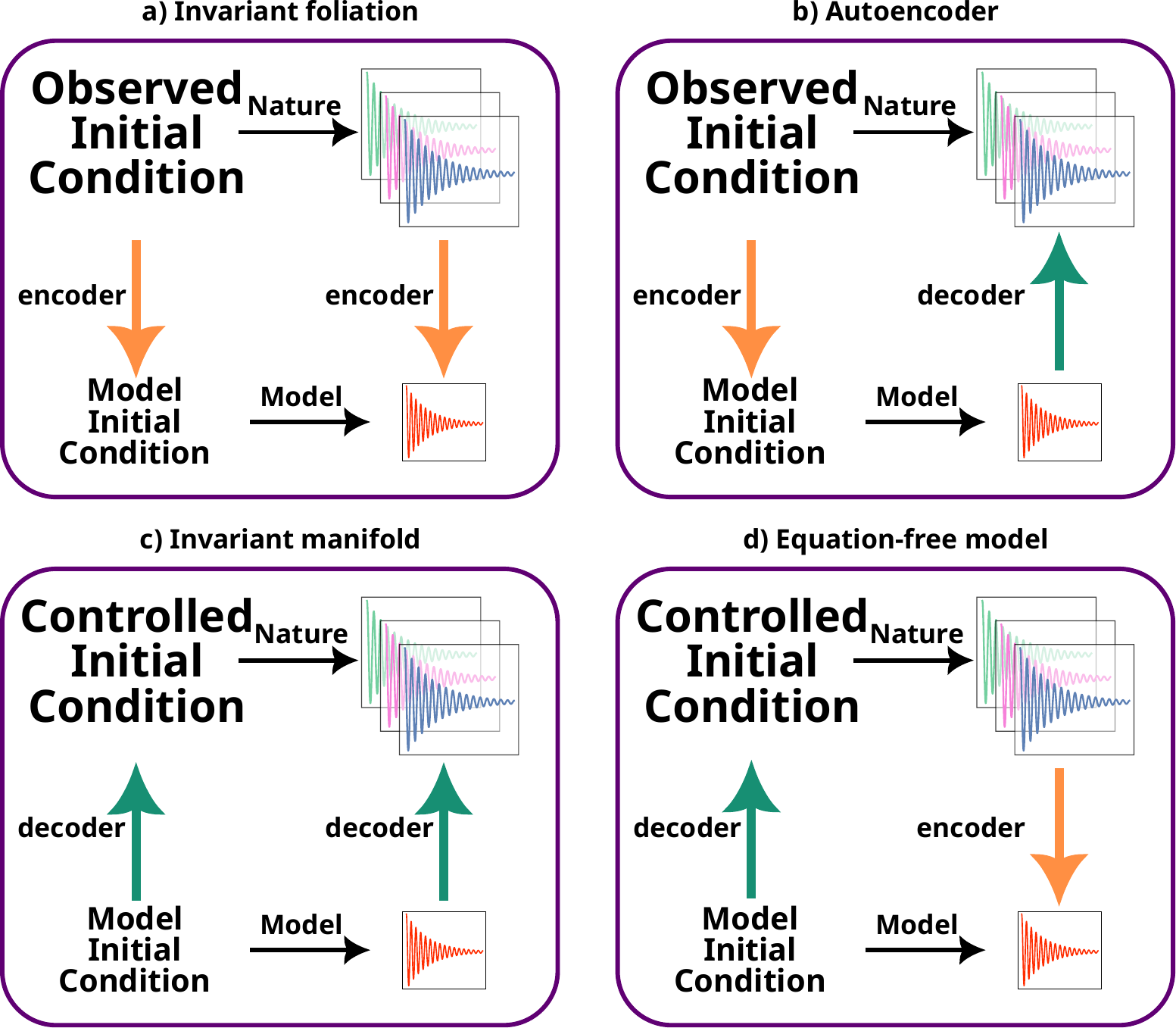}
\par\end{centering}
\caption{\label{fig:four-architectures}Four basic architectures for model
order reduction. a) The observed initial condition and the system
output are both encoded by the same encoder, and the model accuracy
is tested in the latent space; b) the observed initial condition is
encoded but the model prediction is decoded and the prediction is
directly compared to system output, which makes it a generative architecture.
c) The physical initial condition is controlled through a decoder
and the system output is compared to the decoded model output; d)
the physical initial condition is controlled, but model prediction
is compared to the encoded system output.}
\end{figure}

The paper deals with deterministic systems, for which an initial condition
fully determines how the system evolves forward in time. If we think
about how to relate a physical initial condition to model initial
condition, and system output to model prediction, we find four basic
architectures, as depicted in figure \ref{fig:four-architectures}.
The initial condition of the physical system is either observed or
controlled. When the output is observed it needs to be encoded to
become the initial condition of the model. The output of the system
is either encoded and compared with model prediction in the latent
space; or the model generates an object that is directly comparable
to the output of the physical system. This latter approach is used
in generative artificial intelligence. These choices produce four
combinations: invariant foliations \cite{Lawson1974,Roberts89,BatesFoliations2000,AulbachFoliation2003,Szalai2020ISF,Szalai2023Fol},
autoencoders \cite{Kramer1991autoencoder,Champion2019Autoencoder,KaliaMeijerBrunton2021,Cenedese2022NatComm},
invariant manifolds \cite{Fenichel,Bates1998,wiggins2013normally,CabreLlave2003,eldering2013normally,Haro2016},
and equation-free models \cite{Kevrekidis2003,Samey}. 

Here we focus on observed data, which means that at the time of data
analysis, there is no possibility to influence the physical system.
Hence, we are left with either using autoencoders or invariant foliations.
It is shown in \cite{Szalai2023Fol} that autoencoders do not produce
meaningful models, because they only identify where the data is clustered,
which is almost never an invariant manifold. Furthermore, invariant
foliations are unique if their functional representation is sufficiently
smooth and some finite number of resonances are avoided \cite{Szalai2020ISF}.

Invariant foliations were introduced for data-driven modelling in
\cite{Szalai2020ISF} and later proposed in \cite{lecun_path_nodate}
under the name of Joint Encoding Predictive Architecture (JEPA) for
autonomous machine intelligence. It was shown in subsequent papers
\cite{assran2023self,assran2025vjepa2} that JEPA achieves better
accuracy in shorter training time than other tested models. JEPA originates
from siamese networks \cite{Siamese1993,Chicco2021} with the addition
of a predictive model. 

The main motivation for JEPA (or invariant foliations) is to imitate
how humans learn and form concepts from observations. Just like some
aspects of human learning, invariant foliations are not generative.
Since they use a single encoder and a low dimensional model, they
require about half the number of parameters as an autoencoder would
need. Learning occurs in the latent space, which is similar to how
humans conceptualise their learning. For example, a codimension-one
foliation would recognise cars in videos and resolve their speed as
the encoder output. For dynamical systems invariant foliations recognise
modes of vibration or movement and their encoder outputs a signal
that can be independently described by a simple model.

One technical problem with invariant foliations (or JEPA) is that
the encoder can collapse during training if no constraint is set \cite{jing2022collapse}.
In \cite{Szalai2020ISF} a normal form style constraint was used,
which only allowed non-resonant terms in the encoder so that the predictive
model would contain all internal resonances. In \cite{Szalai2023Fol}
a simpler graph-style parametrisation was used, and continued to be
used here.

\subsection{\label{subsec:State-space-reconstruction}State-space reconstruction}

In many cases the system output is lower dimensional than the state
space. In order to reconstruct the state space, delay embedding is
used, inspired by how observability of a system is defined. Locally
and for almost all system outputs it is sufficient to have a delay
length $d$ for which $dm\ge n$, where $n$ is the state space dimension
and $m$ is the dimension of the output signal \cite{Nijmeier1982}.
In the rare case when observability does not hold, delay embedding
forms an invariant foliation \cite{Nijmeier1982}. This provides additional
safety, meaning that the observed dynamics is invariant and meaningful
even if the full phase space cannot be reconstructed. For global observability,
Takens' theorem applies \cite{TakensEmbedding1981}, which states
that if $dm>2n$, almost all system outputs can be used to reconstruct
the state-space globally. As was noted in \cite{Sauer1991}, the required
length of delay for global observability is not necessarily the upper
bound, but falls into the interval $n\le dm\le2n+1$. For dynamics
near a steady state, local observability is sufficient, similar to
how graph style parametrisation is used for invariant manifolds and
foliations. 

Optimal delay embedding can be constructed using principal component
analysis (PCA) \cite{BROOMHEAD1986delayEmbed} or dynamic mode decomposition
(DMD) \cite{SCHMID_2010}. PCA constructs linear combinations of delayed
output values that contain the highest energy de-correlated projections
of the output. DMD on the other hand uses a linear approximation of
the system and calculates its invariant subspaces, where the subspaces
can be chosen based on how much energy the data has when projected
to that given subspace or what subspaces have the lowest residual
\cite{colbrook2023residual}. The number of resolved PCA or DMD modes
$m$ must satisfy $m\ge n$ for local observability, but need not
be more than or $2n+1\le m$ for global observability.

\subsection{PCA and DMD}

It is informative to look into the differences between PCA and DMD
and their relationships to autoencoders and invariant foliations,
respectively. These are linear techniques but they can be used on
nonlinearly embedded data and therefore still able to identify nonlinear
structures. The autoencoder was first introduced as a nonlinear version
of PCA \cite{Kramer1991autoencoder}. Dynamic mode decomposition \cite{SCHMID_2010}
and its variants \cite{Korda_2017} are the linear restrictions of
invariant foliations, because they calculate a linear encoder (in
the form of left eigenvectors) and the identified model is assumed
to be linear. 

The standard assumption for these methods, with the exception of optDMD
\cite{optDMD2018}, is that the data comes as a pair of matrices $\boldsymbol{X},\boldsymbol{Y}\in\mathbb{R}^{d_{X}\times N}$
and the dynamics maps each column of matrix $\boldsymbol{X}$ to the
corresponding column of matrix $\boldsymbol{Y}$. In case of PCA,
a correlation matrix $\boldsymbol{C}=\boldsymbol{X}\boldsymbol{X}^{T}$
is formed and then its singular value decomposition (SVD) is calculated,
such that $\boldsymbol{C}=\boldsymbol{U}\boldsymbol{\Sigma}^{2}\boldsymbol{U}^{T}$.
Matrix $\boldsymbol{\Sigma}$ is diagonal with non-negative entries,
and matrix $\boldsymbol{U}$ is unitary, which consists of a set of
column vectors $\boldsymbol{U}=\begin{pmatrix}\boldsymbol{u}_{1} & \cdots & \boldsymbol{u}_{d_{X}}\end{pmatrix}$.
This identifies correlated patterns in space, along the row indices
of $\boldsymbol{X}$, and weights the most frequently occurring patterns
highly. The data is then reduced to only include the largest $m_{SVD}<d_{X}$
singular values present in the diagonal of matrix $\boldsymbol{\Sigma}$,
which leads to the best (in Frobenius norm) rank-$m_{SVD}$ reduction
of matrix $\boldsymbol{X}$. The decoder and encoder therefore become
$\boldsymbol{U}_{r}=\begin{pmatrix}\boldsymbol{u}_{1} & \cdots & \boldsymbol{u}_{m_{SVD}}\end{pmatrix}$
and $\boldsymbol{U}_{r}^{T}$, respectively, which corresponds to
the architecture in figure \ref{fig:four-architectures}(b). If we
assume a linear system, the reduced order model is found by solving
the least-squares problem
\begin{equation}
\boldsymbol{A}_{r}=\mathop{\mathrm{argmin}}_{\boldsymbol{A}}\left\Vert \boldsymbol{A}\boldsymbol{U}_{r}^{T}\boldsymbol{X}-\boldsymbol{U}_{r}^{T}\boldsymbol{Y}\right\Vert ^{2}.\label{eq:PCA-ROM}
\end{equation}

Dynamic mode decomposition identifies patterns in time and does not
pick up on how much data represents a given temporal pattern. Accuracy,
however depends on the amount of data and can be quantified using
resDMD \cite{colbrook2023residual}. The simplest variant identifies
a linear model
\[
\boldsymbol{A}_{0}=\boldsymbol{Y}\boldsymbol{X}^{T}\left(\boldsymbol{X}\boldsymbol{X}^{T}\right)^{+},
\]
where $^{+}$ stands for the Moore-Penrose pseudo inverse. The linear
model $\boldsymbol{A}_{0}$ (which is assumed to be semisimple) is
then decomposed into a Jordan normal form
\[
\boldsymbol{A}_{0}=\boldsymbol{V}^{-1}\boldsymbol{\Lambda}\boldsymbol{V},
\]
where $\boldsymbol{\Lambda}$ contains the eigenvalues of $\boldsymbol{A}_{0}$
in its diagonal and the rows of $\boldsymbol{V}$ are the corresponding
left eigenvectors of $\boldsymbol{A}_{0}$. The eigenvalues along
the diagonal of $\boldsymbol{\Lambda}$ can be ordered in various
ways: in increasing order of the residual \cite{colbrook2023residual},
decreasing order of the magnitude (slowest dynamics first) or the
amount of data (energy) it captures. The reduced model is then chosen
by the first $m_{DMD}<d_{X}$ eigenvalues. The encoder in figure \ref{fig:four-architectures}(a)
becomes $\boldsymbol{V}_{r}=\begin{pmatrix}\boldsymbol{v}_{1} & \cdots & \boldsymbol{v}_{m_{DMD}}\end{pmatrix}^{T}$,
where $\boldsymbol{v}_{i}$ are the left eigenvectors of $\boldsymbol{A}_{0}$.
The reduced model is
\begin{equation}
\boldsymbol{A}_{r}=\mathop{\mathrm{argmin}}_{\boldsymbol{A}}\left\Vert \boldsymbol{A}\boldsymbol{V}_{r}\boldsymbol{X}-\boldsymbol{V}_{r}\boldsymbol{Y}\right\Vert ^{2}.\label{eq:DMD-ROM}
\end{equation}
This is very similar to \eqref{eq:PCA-ROM}, but now the encoder $\boldsymbol{V}_{r}$
is invariant with respect to the linear dynamics, while $\boldsymbol{U}_{r}$
in \eqref{eq:PCA-ROM} is not. Note that invariance also allows for
direct optimisation
\begin{equation}
\boldsymbol{A}_{r},\boldsymbol{V}_{r}=\mathop{\mathrm{argmin}}_{\boldsymbol{A},\boldsymbol{V},\,\text{s.t.}\,\boldsymbol{V}\boldsymbol{V}^{*}=I}\left\Vert \boldsymbol{A}\boldsymbol{V}\boldsymbol{X}-\boldsymbol{V}\boldsymbol{Y}\right\Vert ^{2}\label{eq:DMD-direct}
\end{equation}
if we restrict the encoder to orthogonal matrices. Due to the orthogonality
of $\boldsymbol{V}$, the resulting $\boldsymbol{A}_{r}$ is no longer
diagonal, but can be of upper Hessenberg form and keeps the same eigenvalues
as $\boldsymbol{A}_{r}$ given by \eqref{eq:DMD-ROM}. This direct
optimisation is not possible for PCA and the identified $\boldsymbol{V}_{r}$
matrix, even if orthogonal, $\boldsymbol{V}_{r}$ is not the same
as $\boldsymbol{U}_{r}^{T}$ in \eqref{eq:PCA-ROM}.

\section{Assumptions and theory}

We assume that the dynamical system we wish to identify is in the
skew-product form 
\begin{equation}
\begin{array}{rl}
\boldsymbol{x}_{k+1} & \negthickspace\negthickspace\negthickspace=\boldsymbol{f}\left(\boldsymbol{x},\boldsymbol{\theta}\right)\\
\boldsymbol{\theta}_{k+1} & \negthickspace\negthickspace\negthickspace=\boldsymbol{g}\left(\boldsymbol{\theta}\right)
\end{array}.\label{eq:Skew_Product}
\end{equation}
where $\boldsymbol{f}\in X\times Y\to X$, $\boldsymbol{g}:Y\to Y$
are real analytic functions defined on the $d_{X}$-dimensional linear
space $X$ and $d_{Y}$-dimensional differentiable manifold $Y$.
The forcing function $\boldsymbol{g}$ is known, while function $\boldsymbol{f}$
is unknown. We further assume that $\boldsymbol{g}$ is a volume preserving
map. Hamiltonian systems generate volume preserving maps according
to Liouville's theorem \cite{RevModPhys.64.795}. By Poincare's recurrence
theorem \cite{Oxtoby1980}, the dynamics of the $\boldsymbol{\theta}$
variable is recurrent, and the Koopman operator associated with $\boldsymbol{g}$
is unitary \cite{colbrook2023mpedmd}. This means that later on we
are able to approximate the forcing with a finite dimensional unitary
transformation (i.e., a unitary matrix).

We assume the existence of an invariant set, called stationary state,
\begin{equation}
\mathcal{T}=\left\{ \boldsymbol{s}\left(\boldsymbol{\theta}\right):\boldsymbol{\theta}\in Y\right\} ,\label{eq:Stationary_State}
\end{equation}
where $\boldsymbol{s}:Y\to X$ satisfies the invariance equation
\begin{equation}
\boldsymbol{s}\left(\boldsymbol{g}\left(\boldsymbol{\theta}\right)\right)=\boldsymbol{f}\left(\boldsymbol{s}\left(\boldsymbol{\theta}\right),\boldsymbol{\theta}\right).\label{eq:S_State_Invariance}
\end{equation}
The set $\mathcal{T}$ can be a periodic orbit, an invariant torus
or a set of fixed points if $Y$ is a parameter space and $\boldsymbol{g}\left(\boldsymbol{\theta}\right)=\boldsymbol{\theta}$.

The data is a set of $N$ trajectories, each of which are $\ell_{j}-\ell_{j-1}$
points long for $j=1,\ldots,N$, where $\ell_{0}=0$. In formulae
the trajectories are 
\begin{align*}
\left(\boldsymbol{x}_{1},\boldsymbol{\theta}_{1}\right),\left(\boldsymbol{x}_{2},\boldsymbol{\theta}_{2}\right),\ldots,\left(\boldsymbol{x}_{\ell_{1}},\boldsymbol{\theta}_{\ell_{1}}\right) & \in X\times Y,\\
 & \vdots\\
\left(\boldsymbol{x}_{\ell_{N-1}+1},\boldsymbol{\theta}_{\ell_{N-1}+1}\right),\left(\boldsymbol{x}_{\ell_{N-1}+2},\boldsymbol{\theta}_{\ell_{N-1}+2}\right),\ldots,\left(\boldsymbol{x}_{\ell_{N}},\boldsymbol{\theta}_{\ell_{N}}\right) & \in X\times Y.
\end{align*}
The data points $\boldsymbol{x}_{k},\boldsymbol{\theta}_{k}$, $k=1,\ldots,\ell_{N}$
representing the deterministic dynamics may be reconstructed from
a system output using delay embedding \cite{TakensEmbedding1981}
as described in section \ref{subsec:State-space-reconstruction}.

In what follows, we transform the system \eqref{eq:Skew_Product}
such that the new state variable becomes 
\begin{equation}
\hat{\boldsymbol{x}}=\boldsymbol{x}-\boldsymbol{s}\left(\boldsymbol{\theta}\right)\label{eq:Data_Translation}
\end{equation}
and then drop the hat. In this translated coordinate system $\boldsymbol{s}\left(\boldsymbol{\theta}\right)=\boldsymbol{0}$,
which simplifies the notation below. This also means that in the transformed
coordinates the stationary state is $\mathcal{T}=\left\{ \boldsymbol{0}:\boldsymbol{\theta}\in Y\right\} $.

\subsection{Invariant foliations}

An invariant foliation can be defined by an encoder (manifold submersion)
$\boldsymbol{u}:X\times Y\to Z$ \cite{Lawson1974}, such that the
leaves of the foliation are
\[
\mathcal{L}\left(\boldsymbol{z},\boldsymbol{\theta}\right)=\left\{ \boldsymbol{x}\in X:\boldsymbol{u}\left(\boldsymbol{x},\boldsymbol{\theta}\right)=\boldsymbol{z}\right\} .
\]
The foliation is the collection of leaves
\[
\mathcal{F}=\left\{ \mathcal{L}\left(\boldsymbol{z},\boldsymbol{\theta}\right):\left(\boldsymbol{z},\boldsymbol{\theta}\right)\in Z\times Y\right\} .
\]
To ensure that the leaves do not intersect each other, the Jacobian
$D_{1}\boldsymbol{u}$ must have full rank everywhere in $X\times Y$.
The space $Z$ is called the latent space, and its dimensionality
is denoted by $d_{Z}$. Foliations are also said to have a co-dimension,
which in our setting is the same as $d_{Z}$. A foliation is invariant
if each leaf is mapped into another leaf by $\boldsymbol{f}$ and
$\boldsymbol{g}$. This is described by the invariance equation
\begin{equation}
\boldsymbol{r}\left(\boldsymbol{u}\left(\boldsymbol{x},\boldsymbol{\theta}\right),\boldsymbol{\theta}\right)=\boldsymbol{u}\left(\boldsymbol{f}\left(\boldsymbol{x},\boldsymbol{\theta}\right),\boldsymbol{g}\left(\boldsymbol{\theta}\right)\right),\label{eq:Invariance_Equation}
\end{equation}
where function $\boldsymbol{r}:Z\times Y\to Z$ is analytic and called
the \emph{conjugate map}. To normalise the foliation, i.e., the representation
of $\boldsymbol{u}$ and $\boldsymbol{r}$ relative to each other,
we prescribe that 
\begin{equation}
\boldsymbol{r}\left(\boldsymbol{0},\boldsymbol{\theta}\right)=\boldsymbol{0}\label{eq:Homogeneous_ROM}
\end{equation}
for all $\boldsymbol{\theta}\in Y$. Given that the stationary state
$\mathcal{T}$ is at the origin, due to the translation \eqref{eq:Data_Translation},
it follows from \eqref{eq:Homogeneous_ROM} that $\boldsymbol{u}\left(\boldsymbol{0},\boldsymbol{\theta}\right)=\boldsymbol{0}$.
Later on, we apply further constraints on the encoder $\boldsymbol{u}$
to prevent it becoming the trivial solution $\boldsymbol{u}\equiv\boldsymbol{0}$
during the solution process of \eqref{eq:Invariance_Equation}. This
is called collapse, and can occur in a variety of ways in self-supervised
learning \cite{jing2022collapse}. For DMD, the collapse is avoided
by making the encoder matrix orthogonal in equation \eqref{eq:DMD-direct}.

\subsection{\label{subsec:Foliations-to-Manifolds}From foliations to manifolds}

Invariant foliations can be used to recover invariant manifolds. We
say that manifold $\mathcal{M}$ is invariant if its forward image
is contained within the manifold. Therefore the leaf of an invariant
foliation that is mapped into itself is an invariant manifold. Due
to the normalising condition \eqref{eq:Homogeneous_ROM}, the leaf
going through the origin $\mathcal{M}=\mathcal{L}\left(\boldsymbol{0},\boldsymbol{\theta}\right)$
is an invariant manifold.

In order to recover an invariant manifold and the dynamics on the
manifold at least two foliations are needed. Let us assume a set of
foliations $\leftidx{^{i\!}}{\mathcal{F}}$, $i=1,\ldots,m_{\mathcal{F}}$
with latent spaces $\leftidx{^{i\!}}{Z}$, conjugate maps $\leftidx{^{i\!}}{\boldsymbol{r}}:\leftidx{^{i\!}}{Z}\times Y\to\leftidx{^{i\!}}{Z}$,
and encoders $\leftidx{^{i\!}}{\boldsymbol{u}}:X\times Y\to\leftidx{^{i\!}}{Z}$.
We further assume that $\sum_{i=1}^{m_{\mathcal{F}}}\dim\leftidx{^{i\!}}{Z}=d_{X}$
and that the linear map
\begin{equation}
\boldsymbol{Q}\left(\boldsymbol{\theta}\right)=\begin{pmatrix}D_{1}\leftidx{^{1\!}}{\boldsymbol{u}}\left(\boldsymbol{0},\boldsymbol{\theta}\right)\\
\vdots\\
D_{1}\leftidx{^{m_{\mathcal{F}}\!}}{\boldsymbol{u}}\left(\boldsymbol{0},\boldsymbol{\theta}\right)
\end{pmatrix}\label{eq:Q-matrix}
\end{equation}
has full rank. Under these assumptions and for an index set 
\[
\mathcal{I}=\left\{ i_{1},\ldots,i_{p}\right\} \subset\left\{ 1,\ldots,m_{\mathcal{F}}\right\} 
\]
of length $p$, there exists an invariant manifold in a neighbourhood
of the stationary state $\mathcal{T}$ defined by 
\[
\leftidx{^{\mathcal{I}\!}}{\mathcal{M}}=\left\{ \left(\boldsymbol{x},\boldsymbol{\theta}\right)\in X\times Y:\leftidx{^{i\!}}{\boldsymbol{u}}\left(\boldsymbol{x},\boldsymbol{\theta}\right)=\boldsymbol{0},\forall i\notin\mathcal{I}\right\} .
\]
The conjugate dynamics on $\leftidx{^{\mathcal{I}\!}}{\mathcal{M}}$
is then given by the map
\[
\boldsymbol{r}\left(\leftidx{^{i_{1}\!}}{\boldsymbol{z}},\ldots,\leftidx{^{i_{p}\!}}{\boldsymbol{z}},\boldsymbol{\theta}\right)=\begin{pmatrix}\leftidx{^{i_{1}\!}}{\boldsymbol{r}}\left(\leftidx{^{i_{1}\!}}{\boldsymbol{z}},\boldsymbol{\theta}\right)\\
\vdots\\
\leftidx{^{i_{p}\!}}{\boldsymbol{r}}\left(\leftidx{^{i_{p}\!}}{\boldsymbol{z}},\boldsymbol{\theta}\right)
\end{pmatrix}.
\]
The invariant manifold has a parametrisation that is compatible with
the parametrisation of the set of foliations, which is given by the
function $\leftidx{^{\mathcal{I}\!}}{\boldsymbol{w}}$ and implicitly
defined by
\[
\leftidx{^{j\!}}{\boldsymbol{u}}\left(\leftidx{^{\mathcal{I}\!}}{\boldsymbol{w}}\left(\leftidx{^{i_{1}\!}}{\boldsymbol{z}},\ldots,\leftidx{^{i_{p}\!}}{\boldsymbol{z}},\boldsymbol{\theta}\right),\boldsymbol{\theta}\right)=\begin{cases}
\leftidx{^{j\!}}{\boldsymbol{z}} & \text{if}\;j\in\mathcal{I}\\
\boldsymbol{0} & \text{if}\;j\notin\mathcal{I}
\end{cases}.
\]
With this definition the manifold can be explicitly written as
\[
\leftidx{^{\mathcal{I}\!}}{\mathcal{M}}=\left\{ \leftidx{^{\mathcal{I}\!}}{\boldsymbol{w}}\left(\leftidx{^{i_{1}\!}}{\boldsymbol{z}},\ldots,\leftidx{^{i_{p}\!}}{\boldsymbol{z}},\boldsymbol{\theta}\right):\leftidx{^{i_{1}\!}}{\boldsymbol{z}},\ldots,\leftidx{^{i_{p}\!}}{\boldsymbol{z}}\in\leftidx{^{i_{1}\!}}{Z}\times\cdots\times\leftidx{^{i_{p}\!}}{Z},\boldsymbol{\theta}\in Y\right\} .
\]

\subsection{\label{subsec:Normal-form}Two-dimensional invariant manifolds and
their normal form}

Assume that there are two foliations, represented by $\leftidx{^{1\!}}{\boldsymbol{r}},$$\leftidx{^{1\!}}{\boldsymbol{u}}$
and $\leftidx{^{2\!}}{\boldsymbol{r}},$$\leftidx{^{2\!}}{\boldsymbol{u}}$
for which matrix $\boldsymbol{Q}$ as given by \eqref{eq:Q-matrix}
has full rank. The first foliation is codimension-two and the Jacobian
$D\leftidx{^{1\!}}{\boldsymbol{r}}\left(\boldsymbol{0},\boldsymbol{\theta}\right)$
is reducible to a matrix with a complex conjugate pair of eigenvalues
as per definition \ref{def:linear-reducible}. Under these assumptions
it is possible to reduce the conjugate dynamics $\leftidx{^{1\!}}{\boldsymbol{r}},$
to
\begin{align}
\rho_{k+1} & =R\left(\rho_{k}\right)\nonumber \\
\beta_{k+1} & =T\left(\rho_{k}\right)\label{eq:Polar-Model-Map}\\
\boldsymbol{\theta}_{k+1} & =\boldsymbol{g}\left(\boldsymbol{\theta}_{k}\right)\nonumber 
\end{align}
in a neighbourhood of the steady state $\mathcal{T}$. Function $R$
describes the damping in the system and function $T$ corresponds
to the instantaneous frequency. These quantities make physical sense
if $\rho$ is proportional to the vibration amplitude and $\beta$
represents an absolute phase in space $X$. Our definition of the
instantaneous relative damping and instantaneous angular frequency
are
\begin{align}
\zeta\left(\rho\right) & =-\frac{\log\frac{\mathrm{d}}{\mathrm{d}\rho}R\left(\rho\right)}{T\left(\rho\right)}\;\;\text{and}\;\;\label{eq:INST-damp}\\
\omega\left(\rho\right) & =\frac{T\left(\rho\right)}{\Delta t},\label{eq:INST-freq}
\end{align}
respectively, where $\Delta t$ is the sampling period of the data.
Elsewhere in the literature \cite{JinBrake2020FreqDamp} the instantaneous
relative damping is defined as 
\begin{equation}
\zeta\left(\rho\right)=-\frac{\log\left(R\left(\rho\right)/\rho\right)}{T\left(\rho\right)}.\label{eq:INST-damp-wrong}
\end{equation}
The two definitions agree at $\rho=0$. We prefer \eqref{eq:INST-damp}
is because it is a truly instantaneous quantity. The alternative definition
\eqref{eq:INST-damp-wrong} is not instantaneous, its value depends
on how the vibration decays in the future.

We now describe how to extract these quantities from the pair of invariant
foliations we assumed to have. The invariance equation of the invariant
manifolds in polar parametrisation is
\begin{equation}
\boldsymbol{w}\left(R\left(\rho\right),T\left(\rho\right),\boldsymbol{g}\left(\boldsymbol{\theta}\right)\right)=\boldsymbol{f}\left(\boldsymbol{w}\left(\rho,\beta,\boldsymbol{\theta}\right),\boldsymbol{\theta}\right).\label{eq:MANIF-polar-inv}
\end{equation}
We do not know $\boldsymbol{f}$, but we can replace it with an invariant
foliation. Applying the encoder $\leftidx{^{1\!}}{\boldsymbol{u}}$
to both sides of \eqref{eq:MANIF-polar-inv} gives us
\begin{align*}
\leftidx{^{1\!}}{\boldsymbol{u}}\left(\boldsymbol{w}\left(R\left(r\right),T\left(r\right),\boldsymbol{g}\left(\boldsymbol{\theta}\right)\right)\right) & =\leftidx{^{1\!}}{\boldsymbol{u}}\left(\boldsymbol{f}\left(\boldsymbol{w}\left(r,\beta,\boldsymbol{\theta}\right),\boldsymbol{\theta}\right)\right).
\end{align*}
Using the foliation invariance equation \eqref{eq:Invariance_Equation},
we get our first equation to solve 
\begin{equation}
\leftidx{^{1\!}}{\boldsymbol{u}}\left(\boldsymbol{w}\left(R\left(r\right),T\left(r\right),\boldsymbol{g}\left(\boldsymbol{\theta}\right)\right)\right)=\leftidx{^{1\!}}{\boldsymbol{r}}\left(\leftidx{^{1\!}}{\boldsymbol{u}}\left(\boldsymbol{w}\left(r,\beta,\boldsymbol{\theta}\right),\boldsymbol{\theta}\right),\boldsymbol{\theta}\right).\label{eq:MANIF-tangential}
\end{equation}
In addition, the invariant manifold is the zero-level set of the second
foliation, hence the equation 
\begin{equation}
\leftidx{^{2\!}}{\boldsymbol{u}}\left(\boldsymbol{w}\left(r,\beta,\boldsymbol{\theta}\right)\right)=\boldsymbol{0}\label{eq:MANIF-normal}
\end{equation}
also needs to hold. Now we need to ensure that $\rho$ is proportional
to the vibration amplitude, hence we impose the constraint 
\begin{equation}
\int_{Y}\int_{0}^{2\pi}\left\langle D_{1}\boldsymbol{w}\left(\rho,\beta,\boldsymbol{\theta}\right),\boldsymbol{w}\left(\rho,\beta,\boldsymbol{\theta}\right)-\boldsymbol{w}\left(0,\beta,\boldsymbol{\theta}\right)\right\rangle \mathrm{d}\beta\mathrm{d}\boldsymbol{\theta}=\rho.\label{eq:MANIF-ampl}
\end{equation}
To make sure that there is no phase distortion in space $X$ with
respect to the $\beta$ variable, we also need to impose
\begin{equation}
\int_{Y}\int_{0}^{2\pi}\left\langle D_{1}\boldsymbol{w}\left(r,\beta,\boldsymbol{\theta}\right),D_{2}\boldsymbol{w}\left(r,\beta,\boldsymbol{\theta}\right)\right\rangle \mathrm{d}\beta\mathrm{d}\boldsymbol{\theta}=0.\label{eq:MANIF-phase}
\end{equation}
For a detailed derivation of the conditions \eqref{eq:MANIF-ampl}
and \eqref{eq:MANIF-phase}, see \cite{Szalai2023Fol}. Equations
\eqref{eq:MANIF-tangential}, \eqref{eq:MANIF-normal}, \eqref{eq:MANIF-ampl}
and \eqref{eq:MANIF-phase} form a system that uniquely determines
functions $R$,\textbf{ }$T$ and $\boldsymbol{w}$. In this paper
we discretise the equations through a collocation technique and obtain
the solution using Newton's method.
\begin{rem}
\label{rem:MAP-backbones}If we know map $\boldsymbol{f}$ then the
three equations \eqref{eq:MANIF-polar-inv}, \eqref{eq:MANIF-ampl}
and \eqref{eq:MANIF-phase} can also be used to calculate the two-dimensional
invariant manifold and the instantaneous frequency and damping ratio
of the dynamics on the manifold.
\end{rem}

\begin{rem}
\label{rem:ODE-backbones}In case of a differential equation
\begin{align*}
\dot{\boldsymbol{x}} & =\boldsymbol{F}\left(\boldsymbol{x},\boldsymbol{\theta}\right)\\
\dot{\boldsymbol{\theta}} & =\boldsymbol{G}\left(\boldsymbol{\theta}\right)
\end{align*}
we can reduce the dynamics to 
\begin{align}
\dot{\rho} & =\hat{R}\left(\rho\right)\nonumber \\
\dot{\beta} & =\hat{T}\left(\rho\right)\label{eq:Polar-Model-ODE}\\
\dot{\boldsymbol{\theta}} & =\boldsymbol{G}\left(\boldsymbol{\theta}\right)\nonumber 
\end{align}
for which the invariance equation is
\[
D_{1}\boldsymbol{w}\left(\rho,\beta,\boldsymbol{\theta}\right)\hat{R}\left(\rho\right)+D_{2}\boldsymbol{w}\left(\rho,\beta,\boldsymbol{\theta}\right)\hat{T}\left(\rho\right)+D_{3}\boldsymbol{w}\left(\rho,\beta,\boldsymbol{\theta}\right)\boldsymbol{G}\left(\boldsymbol{\theta}\right)=\boldsymbol{F}\left(\boldsymbol{w}\left(\rho,\beta,\boldsymbol{\theta}\right),\boldsymbol{\theta}\right).
\]
The two constraints \eqref{eq:MANIF-ampl} and \eqref{eq:MANIF-phase}
remain the same, while the instantaneous relative damping and instantaneous
angular frequency become
\begin{align}
\zeta\left(\rho\right) & =-\frac{\frac{\mathrm{d}}{\mathrm{d}\rho}R\left(\rho\right)}{T\left(\rho\right)}\;\;\text{and}\;\;\label{eq:INST-damp-ode}\\
\omega\left(\rho\right) & =T\left(\rho\right),\label{eq:INST-freq-ode}
\end{align}
respectively.
\end{rem}

\subsection{\label{subsec:Fitting-a-foliation}Fitting a foliation to trajectories}

For two consecutive data points, the invariance equation \eqref{eq:Invariance_Equation}
reads
\[
\boldsymbol{r}\left(\boldsymbol{u}\left(\boldsymbol{x}_{k},\boldsymbol{\theta}_{k}\right),\boldsymbol{\theta}_{k}\right)=\boldsymbol{u}\left(\boldsymbol{x}_{k+1},\boldsymbol{\theta}_{k+1}\right)
\]
when $k\neq\ell_{j},$$j=0,\ldots,N-1$. However, data comes as a
set of trajectories and the invariance equation \eqref{eq:Invariance_Equation}
must apply over a full trajectory. This means that trajectories in
the latent space $Z$ produced by map $\boldsymbol{r}$ must match
up with the data encoded through function $\boldsymbol{u}$. In particular,
the model in the latent space is
\[
\begin{array}{rl}
\boldsymbol{z}_{k+1} & \negthickspace\negthickspace\negthickspace=\boldsymbol{r}\left(\boldsymbol{z}_{k},\boldsymbol{\theta}_{k}\right)\\
\boldsymbol{\theta}_{k+1} & \negthickspace\negthickspace\negthickspace=\boldsymbol{g}\left(\boldsymbol{\theta}\right)
\end{array},
\]
which can be extended to describe longer trajectories 
\[
\boldsymbol{z}_{k+l}=\boldsymbol{r}\left(\cdots\boldsymbol{r}\left(\boldsymbol{z}_{k},\boldsymbol{\theta}_{k}\right)\cdots,\boldsymbol{g}^{l-1}\left(\boldsymbol{\theta}_{k}\right)\right)\stackrel{\mathrm{def}}{=}\boldsymbol{r}^{l}\left(\boldsymbol{z}_{k},\boldsymbol{\theta}_{k}\right).
\]
Invariance along a trajectory is now written as
\begin{equation}
\boldsymbol{r}^{l-1}\left(\bar{\boldsymbol{z}}_{j+1},\boldsymbol{\theta}_{\ell_{j}}\right)=\boldsymbol{u}\left(\boldsymbol{x}_{\ell_{j}+l},\boldsymbol{\theta}_{\ell_{j}+l}\right),\qquad l=1,\ldots,\ell_{j+1}-\ell_{j},j=0,\ldots,N-1,\label{eq:Multistep_Invariance}
\end{equation}
where 
\begin{equation}
\bar{\boldsymbol{z}}_{j+1}=\boldsymbol{u}\left(\boldsymbol{x}_{\ell_{j}+1},\boldsymbol{\theta}_{\ell_{j}+1}\right)\label{eq:Latent_IC}
\end{equation}
is the initial condition of the trajectory in the latent space $Z$.
When solving the invariance equation \eqref{eq:Multistep_Invariance}
by optimisation, $\bar{\boldsymbol{z}}_{j+1}$ is kept as an unknown,
which acknowledges the fact that the initial data point $\boldsymbol{x}_{\ell_{j}+1},\boldsymbol{\theta}_{\ell_{j}+1}$
might be inaccurate and should not unduly influence the rest of the
latent trajectory.

\begin{figure}
\begin{centering}
\[\begin{tikzcd}
	& {} \\
	{\boldsymbol{x}_{\ell_j+1},\boldsymbol{\theta}_{\ell_j+1}} && {\boldsymbol{x}_{\ell_j+2},\boldsymbol{\theta}_{\ell_j+2}} && {\cdots } && {\boldsymbol{x}_{\ell_{j+1}},\boldsymbol{\theta}_{\ell_{j+1}}} \\
	\\
	{{\Large\textcolor{blue}{\bar{\boldsymbol{z}}_{j+1}}}} && {\boldsymbol{z}_{\ell_j+2}} && {\cdots } && {\boldsymbol{z}_{\ell_{j+1}}}
	\arrow["{\boldsymbol{f},\boldsymbol{g}}", from=2-1, to=2-3]
	\arrow["{{\Large\textcolor{blue}{\boldsymbol{u}}}}"{description}, from=2-1, to=4-1]
	\arrow["{\boldsymbol{f},\boldsymbol{g}}", from=2-3, to=2-5]
	\arrow["{{\Large\textcolor{blue}{\boldsymbol{u}}}}"{description}, from=2-3, to=4-3]
	\arrow["{\boldsymbol{f},\boldsymbol{g}}", from=2-5, to=2-7]
	\arrow["{{\Large\textcolor{blue}{\boldsymbol{u}}}}"{description}, squiggly, tail, two heads, from=2-5, to=4-5]
	\arrow["{{\Large\textcolor{blue}{\boldsymbol{u}}}}"{description}, from=2-7, to=4-7]
	\arrow["{{\Large\textcolor{blue}{\boldsymbol{r}}}}", from=4-1, to=4-3]
	\arrow["{{\Large\textcolor{blue}{\boldsymbol{r}}}}", from=4-3, to=4-5]
	\arrow["{{\Large\textcolor{blue}{\boldsymbol{r}}}}", from=4-5, to=4-7]
\end{tikzcd}\]
\par\end{centering}
\caption{\label{fig:Trajectory_Invariance}A diagrammatic view of the invariance
equation \eqref{eq:Multistep_Invariance}. All data points along a
trajectory are encoded using function $\boldsymbol{u}$ and the encoded
values are compared to the model trajectory produced by map $\boldsymbol{r}$.
The initial condition of the model trajectory $\bar{\boldsymbol{z}}_{j+1}$
is unknown and so are the parameters of function $\boldsymbol{u}$
and $\boldsymbol{r}$, which are typeset in blue.}

\end{figure}
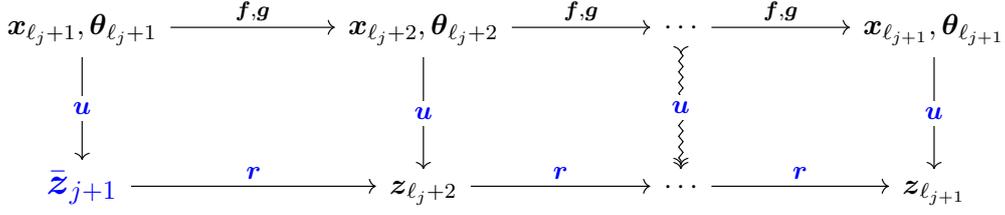

When using data to determine the functions $\boldsymbol{r}$ and $\boldsymbol{u}$,
we minimise the loss function
\begin{equation}
L\left(\boldsymbol{r},\boldsymbol{u},\bar{\boldsymbol{z}}_{1},\ldots,\bar{\boldsymbol{z}}_{N}\right)=\frac{1}{2}\sum_{j=0}^{N-1}\sum_{l=1}^{\ell_{j+1}-\ell_{j}}\sigma_{\epsilon}\left(\left|\boldsymbol{x}_{\ell_{j}+l}-\boldsymbol{s}\left(\boldsymbol{\theta}_{\ell_{j}+l}\right)\right|\right)\left|\boldsymbol{r}^{l-1}\left(\bar{\boldsymbol{z}}_{j},\boldsymbol{\theta}_{\ell_{j}}\right)-\boldsymbol{u}\left(\boldsymbol{x}_{\ell_{j}+l},\boldsymbol{\theta}_{\ell_{j}+l}\right)\right|^{2},\label{eq:Abstract_Loss_Function}
\end{equation}
where $\sigma_{\epsilon}$ is a strictly positive weight function. 

When identifying the invariant foliation, we also search for the best
initial conditions $\bar{\boldsymbol{z}}_{1},\ldots,\bar{\boldsymbol{z}}_{N}$
that create the best fitting trajectory instead of prescribing \eqref{eq:Latent_IC}.
This avoids the large influence that any possible inaccuracy of the
first point of an observed trajectory$\boldsymbol{x}_{\ell_{j}+1},\boldsymbol{\theta}_{\ell_{j}+1}$
has on the latent model trajectory. The weight function can be used
to assign higher importance to some data points than others based
on their amplitude about the steady state.

\begin{figure}
\begin{centering}
\includegraphics[width=0.49\textwidth]{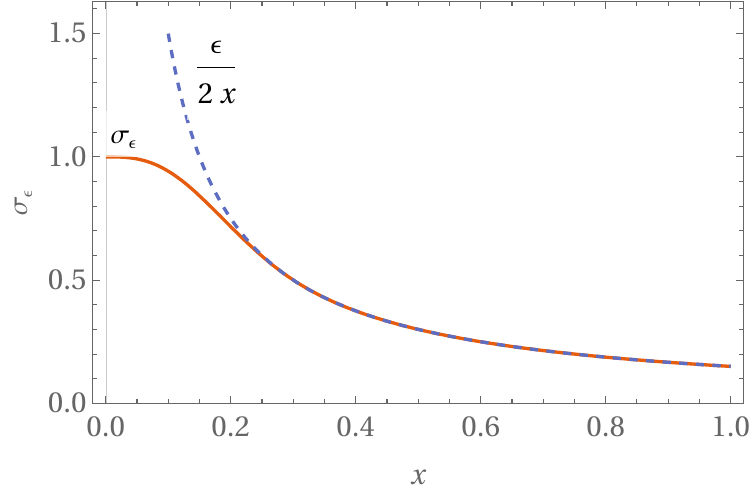}
\par\end{centering}
\caption{\label{fig:ScalingFunction}The graph of the weight function \eqref{eq:ScalingFunction}
for $\epsilon=0.2$.}
\end{figure}

Our pick of the weight function approximates the inverse of the absolute
distance from the steady state for distances greater than $\epsilon$
and makes sure that it avoids the singularity as the distance tends
to zero. The scaling function is
\begin{equation}
\sigma_{\epsilon}\left(x\right)=\begin{cases}
\left(1+\frac{2x^{3}}{\epsilon^{3}}-\frac{x^{4}}{\epsilon^{4}}\right)^{-1} & 0\le x<\epsilon\\
\frac{\epsilon}{2x} & \epsilon\le x
\end{cases},\label{eq:ScalingFunction}
\end{equation}
which is illustrated in figure \ref{fig:ScalingFunction}. The value
of $\epsilon$ is chosen based on the level of noise in the data.
We also define the relative error for a given data point as 
\begin{equation}
E_{\mathit{rel}}\left(\boldsymbol{x}_{\ell_{j}+l},\boldsymbol{\theta}_{\ell_{j}+l}\right)=\frac{2}{\epsilon}\sigma_{\epsilon}\left(\left|\boldsymbol{x}_{\ell_{j}+l}-\boldsymbol{s}\left(\boldsymbol{\theta}_{\ell_{j}+l}\right)\right|\right)\left|\boldsymbol{r}^{l-1}\left(\bar{\boldsymbol{z}}_{j},\boldsymbol{\theta}_{\ell_{j}}\right)-\boldsymbol{u}\left(\boldsymbol{x}_{\ell_{j}+l},\boldsymbol{\theta}_{\ell_{j}+l}\right)\right|.\label{eq:Relative_Error}
\end{equation}
Note that the last factor of \eqref{eq:Relative_Error} is not squared.

\subsection{\label{subsec:Vector-Bundles}Linear decomposition}

To decide what dynamics should the conjugate map $\boldsymbol{r}$
resolve, we investigate the linear dynamics of the system about the
stationary state $\mathcal{T}$. The linearised system about the stationary
state is 
\begin{equation}
\begin{array}{rl}
\boldsymbol{x}_{k+1} & \negthickspace\negthickspace\negthickspace=\boldsymbol{A}\left(\boldsymbol{\theta}_{k}\right)\boldsymbol{x}_{k}\\
\boldsymbol{\theta}_{k+1} & \negthickspace\negthickspace\negthickspace=\boldsymbol{g}\left(\boldsymbol{\theta}_{k}\right)
\end{array},\;k=1,2,\ldots,\label{eq:ED-linsys-1}
\end{equation}
where
\[
\boldsymbol{A}\left(\boldsymbol{\theta}\right)=D_{1}\boldsymbol{f}\left(\boldsymbol{s}\left(\boldsymbol{\theta}\right),\boldsymbol{\theta}\right).
\]
Instead of eigenvectors and eigenvalues, we need to use invariant
vector bundles to decompose the dynamics of the linear system \eqref{eq:ED-linsys-1}.
Invariant foliations require the use of left vector bundles that satisfy
the linear invariance equation 
\begin{equation}
\boldsymbol{\Lambda}_{j}\left(\boldsymbol{\theta}\right)\boldsymbol{U}_{j}\left(\boldsymbol{\theta}\right)=\boldsymbol{U}_{j}\left(\boldsymbol{g}\left(\boldsymbol{\theta}\right)\right)\boldsymbol{A}\left(\boldsymbol{\theta}\right),\label{eq:ED-left-bundle}
\end{equation}
where $\boldsymbol{U}_{j}:Y\to L\left(X,Z_{j}\right)$  and $\boldsymbol{\Lambda}_{j}:Y\to L\left(Z_{j},Z_{j}\right)$
are analytic matrix valued functions and $Z_{j}$ are a low-dimensional
vector spaces, where $j=1,\ldots,m_{\Sigma}$. The decomposition is
complete when $X$ is isomorphic to $\bigoplus_{j=1}^{m_{\Sigma}}Z_{j}$,
and $\bigoplus_{j=1}^{m_{\Sigma}}\left(\ker\boldsymbol{U}_{j}\left(\boldsymbol{\theta}\right)\right)^{\perp}=X$
for all $\boldsymbol{\theta}\in\mathbb{T}^{d_{Y}}$. To characterise
the linear dynamics, we use exponential dichotomies. The matrix $\boldsymbol{\Lambda}_{j}$
has an exponential dichotomy for $\rho\in\mathbb{R}^{+}$ if there
exists $C>0$ such that the inequalities 
\begin{align*}
\left|\boldsymbol{\Lambda}_{j}\left(\boldsymbol{g}^{k-1}\left(\boldsymbol{\theta}\right)\right)\cdots\boldsymbol{\Lambda}_{j}\left(\boldsymbol{g}\left(\boldsymbol{\theta}\right)\right)\boldsymbol{\Lambda}_{j}\left(\boldsymbol{\theta}\right)\right| & \le C\rho^{k}, & \text{and}\\
\left|\boldsymbol{\Lambda}_{j}^{-1}\left(\boldsymbol{g}^{-k}\left(\boldsymbol{\theta}\right)\right)\cdots\boldsymbol{\Lambda}_{j}^{-1}\left(\boldsymbol{g}^{-1}\left(\boldsymbol{\theta}\right)\right)\right| & \le C\rho^{-k}, & \forall k\in\mathbb{N}
\end{align*}
hold. Our spectral decomposition is such that $\boldsymbol{\Lambda}_{j}$
does not have exponential dichotomy when $\rho\in\Sigma_{j}=\left[\alpha_{j},\beta_{j}\right]$
but does anywhere else in $\mathbb{R}^{+}$. The sets $\Sigma_{j}=\left[\alpha_{j},\beta_{j}\right]$
are called spectral intervals, that are pairwise disjoint, i.e. $\Sigma_{j}\cap\Sigma_{k}=\emptyset$
for $j\neq k$. 

For a constant matrix $\boldsymbol{A}$, the spectral intervals reduce
to points ($\alpha_{j}=\beta_{j}$), which are the magnitudes of the
eigenvalues of $\boldsymbol{A}$. In addition, the row vectors of
$\boldsymbol{U}_{j}$ span the same space as the left eigenvectors
of $\boldsymbol{A}$ corresponding to the eigenvalues that have the
magnitude $\alpha_{j}=\beta_{j}$. In section \ref{subsec:Bundles}
we detail how to find $\boldsymbol{\Lambda}_{j}$ and $\boldsymbol{U}_{j}$
numerically by discretising \eqref{eq:ED-left-bundle} and turning
it into an eigenvalue problem.

Note that \eqref{eq:ED-left-bundle} is the linearised version of
the invariance equation \eqref{eq:Invariance_Equation} and therefore
the linear approximation of the foliation can be constructed such
that 
\begin{align}
\leftidx{^{\mathcal{I}\!}}{\boldsymbol{u}}\left(\boldsymbol{z},\boldsymbol{\theta}\right) & =\begin{pmatrix}\boldsymbol{\Lambda}_{i_{1}}\left(\boldsymbol{\theta}\right) &  &  & \boldsymbol{0}\\
 & \boldsymbol{\Lambda}_{i_{2}}\left(\boldsymbol{\theta}\right)\\
 &  & \ddots\\
\boldsymbol{0} &  &  & \boldsymbol{\Lambda}_{i_{p}}\left(\boldsymbol{\theta}\right)
\end{pmatrix}\boldsymbol{z}+\mathcal{O}\left(\left|\boldsymbol{z}\right|^{2}\right),\label{eq:Linear_ROM}\\
\leftidx{^{\mathcal{I}\!}}{\boldsymbol{u}}\left(\boldsymbol{x},\boldsymbol{\theta}\right) & =\begin{pmatrix}\boldsymbol{U}_{i_{1}}\left(\boldsymbol{\theta}\right)\\
\boldsymbol{U}_{i_{2}}\left(\boldsymbol{\theta}\right)\\
\vdots\\
\boldsymbol{U}_{i_{p}}\left(\boldsymbol{\theta}\right)
\end{pmatrix}\boldsymbol{x}+\mathcal{O}\left(\left|\boldsymbol{x}\right|^{2}\right),\label{eq:Linear_Encoder}
\end{align}
where $\mathcal{I}=\left\{ i_{1},\ldots,i_{p}\right\} $ is an index
set and $i_{1},\ldots,i_{p}\in\left\{ 1,2,\ldots,m_{\Sigma}\right\} $
are suitably chosen non-repeating indices. In our examples a single
index $i_{1}$ is chosen such that $\boldsymbol{\Lambda}_{i_{1}}$
correspond to a two-dimensional subspace with oscillatory dynamics,
but the accompanying software \cite{InvarianModels2025} allows for
any index set. When minimising the loss function \eqref{eq:Abstract_Loss_Function},
formulae \eqref{eq:Linear_ROM}, \eqref{eq:Linear_Encoder}, and the
initial conditions \eqref{eq:Latent_IC} serve as the starting point.

\begin{defn}
\label{def:linear-reducible}When $\Lambda_{j}$ can be made independent
of $\boldsymbol{\theta}$, we call the linear system \eqref{eq:ED-linsys-1}
\emph{reducible}. Otherwise \eqref{eq:ED-linsys-1} is \emph{irreducible.}
\end{defn}

\begin{rem}
Note that if $\boldsymbol{g}\left(\boldsymbol{\theta}\right)=\boldsymbol{\theta}$,
the variable $\boldsymbol{\theta}$ is just a set of parameters. In
this case $\boldsymbol{A}\left(\boldsymbol{\theta}\right)$ remains
constant for any given trajectory, though it varies from one trajectory
to another. The spectral intervals are therefore the union of all
eigenvalues of $\boldsymbol{A}\left(\boldsymbol{\theta}\right)$
\[
\Sigma=\bigcup_{j=1}^{m_{\Sigma}}\Sigma_{j}=\left\{ \mathrm{eig}\boldsymbol{A}\left(\boldsymbol{\theta}\right):\boldsymbol{\theta}\in Y\right\} .
\]
\end{rem}

\subsection{Periodic and quasi-periodic forcing}

This section focuses on the case when $Y=\mathbb{T}^{d_{Y}}$ is the
$d_{Y}$-dimensional torus and the forcing is the constant rotation
on the torus, that is 
\begin{equation}
\boldsymbol{g}\left(\boldsymbol{\theta}\right)=\boldsymbol{\theta}+\boldsymbol{\omega}\label{eq:Rigid_Rotation}
\end{equation}
where $\boldsymbol{\omega}\in\mathbb{T}^{d_{Y}}$ is a constant. In
this case $\mathcal{T}$ is also a $d_{Y}$-dimensional torus. The
following theorem sets conditions for the existence and uniqueness
of the invariant foliation.
\begin{thm}
\label{thm:Existence-Uniqueness}Assume an invariant torus \eqref{eq:Stationary_State}
and a linearised system about the torus in the form of \eqref{eq:ED-linsys-1}.
Also assume that the linear system has dichotomy spectral intervals
$\Sigma_{j}=\left[\alpha_{j},\beta_{j}\right]$, $1\le j\le m_{\Sigma}$
and that $\max\beta_{j}<1$. Pick one or more spectral intervals using
a non-empty index set 
\[
\mathcal{I}\subset\left\{ 1,2,\ldots,m_{\Sigma}\right\} 
\]
with the condition that $\alpha_{j}\neq0$ for all $j\in\mathcal{I}$
so that linear system \eqref{eq:ED-linsys-1} restricted to the subset
of the spectrum 
\[
\Sigma_{\mathcal{F}}=\bigcup_{j\in\mathcal{I}}\Sigma_{j}
\]
is invertible. Define the spectral quotient as 
\[
\beth_{\mathcal{I}}=\frac{\min_{j\in\mathcal{I}}\log\alpha_{j}}{\log\beta_{m_{\Sigma}}}.
\]
If the non-resonance conditions 
\begin{equation}
1\notin\left[\beta_{i_{0}}^{-1}\alpha_{i_{1}}\cdots\alpha_{i_{j}},\alpha_{i_{0}}^{-1}\beta_{i_{1}}\cdots\beta_{i_{j}}\right]\label{eq:FOIL-non-resonance-1}
\end{equation}
hold for $i_{0}\in\mathcal{I}$, $i_{1},\ldots,i_{j}\in\left\{ 1,2,\ldots,m_{\Sigma}\right\} $
and $2\le j<\beth_{\mathcal{I}}+1$ with some $i_{k}\notin\mathcal{I}$
then
\begin{enumerate}
\item there exists a unique invariant foliation defined by analytic functions
$\boldsymbol{u}$ and $\boldsymbol{r}$ satisfying equation \eqref{eq:Invariance_Equation}
in a sufficiently small neighbourhood of the torus $\mathcal{T}$
such that equations \eqref{eq:Linear_ROM}, \eqref{eq:Linear_Encoder}
hold;
\item the nonlinear map $\boldsymbol{R}$ is a polynomial, which in its
simplest form contains terms for which the internal non-resonance
conditions \eqref{eq:FOIL-non-resonance-1} with $i_{0},i_{1},\ldots,i_{j}\in\mathcal{I}$
and $2\le j<\beth_{\mathcal{I}}+1$ does not hold.
\end{enumerate}
\end{thm}

\begin{proof}
The details can be found in the paper \cite{SzalaiForcedTheory2024}.
\end{proof}
\begin{rem}
The theorem can be extended to forcing through integrable volume preserving
maps. These maps can always be written in action-angle form, where
$\boldsymbol{\theta}=\left(\boldsymbol{\theta}_{1},\boldsymbol{\theta}_{2}\right)^{T}$such
that
\[
\boldsymbol{g}\left(\boldsymbol{\theta}\right)=\begin{pmatrix}\boldsymbol{\theta}_{1}\\
\boldsymbol{\theta}_{2}+\boldsymbol{\omega}\left(\boldsymbol{\theta}_{1}\right)
\end{pmatrix}.
\]
Given that $\boldsymbol{\theta}_{1}$ is constant and $\boldsymbol{\theta}_{2}$
is a rigid rotation, theorem \ref{thm:Existence-Uniqueness} still
holds.
\end{rem}

\begin{rem}
In many cases an autonomous conjugate map can be chosen, such that
$\boldsymbol{r}\left(\boldsymbol{z},\boldsymbol{\theta}\right)=\boldsymbol{r}\left(\boldsymbol{z}\right)$,
which we call the \emph{nonlinearly reducible} case. Nonlinear reducibility
requires further non-resonance conditions in addition to \eqref{eq:FOIL-non-resonance-1}.
\end{rem}

\begin{rem}
In case of nonlinear reducibility, it is possible to use autonomous
invariant foliations as long as the data contains information about
how the forcing is produced. The forcing dynamics $\boldsymbol{g}$
can then be directly uncovered from data by calculating its own invariant
foliation. This has been tried for some examples (not shown), produces
the same result. The post-processing is more involved and will not
be detailed here.
\end{rem}

\section{Numerical approximation of the foliation}

This section deals with numerical representation of functions $\boldsymbol{r}$
and $\boldsymbol{u}$ and how to calculate their initial linear approximations
\eqref{eq:Linear_ROM}, \eqref{eq:Linear_Encoder}. A suitable choice
of representation avoids both overfitting and underfitting, and depends
on the data. There are many ways to represent a function, using linear
combination of basis functions, constructive approximation (e.g. wavelets)
\cite{devore1998}, neural networks \cite{ElbrachterBolcskei2021}
and even non-parametric regression \cite{NonparametricRegression2008}.
In this section we explore only some limited possibilities.

\paragraph{Notation.}

Entries of a vector $\boldsymbol{v}\in\mathbb{R}^{n}$ are denoted
by $v_{i}$, entries of a matrix $\boldsymbol{A}\in\mathbb{R}^{n\times m}$
are denoted by $A_{ij}$ and entries of a tensor of order three $\boldsymbol{B}\in\mathbb{R}^{n\times m\times p}$
are denoted by $B_{ijk}$, where the indices are $i=1,\ldots,n$,
$j=1,\ldots,m$, and $k=1,\ldots,p$. We use Einstein notation when
it comes to tensor contraction. For example the matrix-vector product
is written as $w_{i}=A_{ij}v_{j}$, where summation is implicitly
understood over index $j$. In some cases summation is not carried
out over an index that occurs in multiple factors, which is signalled
by underlined indices, as in $B_{il}=A_{ijk}v_{j\underline{l}}w_{k\underline{l}}$
index $l$ is not summed over. When we want to indicate that a group
of indices forms rows indices and another group of indices forms columns
indices of an array, we put parentheses around those indices. The
matrix whose row indices are $i,j$ and column indices are $k,l$
is denoted by $A_{\left(ij\right)\left(kl\right)}$.

\paragraph{~}

To approximate the conjugate map $\boldsymbol{r}$, we use a linear
combination of basis functions that are taken from a tensor product
space of functions in the two variables on $X\times Y$. In the first
variable we use a polynomial basis with $n_{Z}$ monomials and in
the second variable we use $n_{Y}$ basis functions that are problem
specific, e.g. the Dirichlet kernel for quasi-periodic focing. Our
approximation is written as
\[
r_{i}\left(\boldsymbol{z},\boldsymbol{\theta}\right)=R_{ijk}\phi_{j}\left(\boldsymbol{z}\right)\psi_{k}\left(\boldsymbol{\theta}\right),
\]
where $R_{ijk}$ are the parameters to be identified. In accordance
with our notation
\[
\boldsymbol{\Phi}\left(\boldsymbol{z}\right)=\begin{pmatrix}\phi_{1}\left(\boldsymbol{z}\right)\\
\phi_{2}\left(\boldsymbol{z}\right)\\
\vdots\\
\phi_{n_{Z}}\left(\boldsymbol{z}\right)
\end{pmatrix},\quad\boldsymbol{\Psi}\left(\boldsymbol{\theta}\right)=\begin{pmatrix}\psi_{1}\left(\boldsymbol{\theta}\right)\\
\psi_{2}\left(\boldsymbol{\theta}\right)\\
\vdots\\
\psi_{n_{Y}}\left(\boldsymbol{\theta}\right)
\end{pmatrix}.
\]

Given that our data is fixed, we can pre-compute 
\[
\boldsymbol{\alpha}_{j}=\boldsymbol{\Psi}\left(\boldsymbol{\theta}_{j}\right),\quad j=1,\ldots,\ell_{N}
\]
to replace the effect of forcing. Moreover, $\boldsymbol{\alpha}_{j}\in\mathbb{R}^{n_{Y}}$
has its own dynamics 
\[
\boldsymbol{\alpha}_{j+1}=\boldsymbol{\Gamma}\left(\boldsymbol{\alpha}_{j}\right),
\]
where $\boldsymbol{\Gamma}$ satisfies the invariance equation
\[
\boldsymbol{\Gamma}\left(\boldsymbol{\Psi}\left(\boldsymbol{\theta}\right)\right)=\boldsymbol{\Phi}\left(\boldsymbol{g}\left(\boldsymbol{\theta}\right)\right).
\]
As in dynamic mode decomposition \cite{SCHMID_2010}, we approximate
$\boldsymbol{\Gamma}$ by a linear function, that is, $\boldsymbol{\Gamma}\left(\boldsymbol{\alpha}\right)=\boldsymbol{\Omega}\boldsymbol{\alpha}$.
This approximation makes sense, because $\boldsymbol{\Omega}$ converges
to the Koopman operator as $n_{Y}\to\infty$ \cite{KlusKoopma2016,Korda_2017}.
Hence we use least squares to fit $\boldsymbol{\Omega}$ to the data,
which yields
\begin{equation}
\boldsymbol{\Omega}=\left(\sum_{j=0}^{N-1}\sum_{k=\ell_{j}}^{\ell_{j+1}-1}\boldsymbol{\alpha}_{k+1}\boldsymbol{\alpha}_{k}^{T}\right)\left(\sum_{j=0}^{N-1}\sum_{k=\ell_{j}}^{\ell_{j+1}-1}\boldsymbol{\alpha}_{k}\boldsymbol{\alpha}_{k}^{T}\right)^{-1}.\label{eq:Omega_LSQ}
\end{equation}
We note that $\boldsymbol{\Omega}$ must be a unitary matrix, because
we assumed that the forcing dynamics is volume preserving. For more
accurate results, measure preserving dynamic mode decomposition can
be used \cite{colbrook2023mpedmd} instead of formula \eqref{eq:Omega_LSQ}.

In case $Y=\mathbb{T}$ we choose a positive odd integer $n_{Y}$
to form a uniform grid
\begin{equation}
\vartheta_{j}=\left(j-1\right)\frac{2\pi}{n_{Y}},\quad j=1,\ldots,n_{Y}\label{eq:circle-colloc}
\end{equation}
and define the library of functions by
\[
\psi_{j}^{n_{Y}}\left(\theta\right)=\frac{1}{n_{Y}}\gamma^{n_{Y}}\left(\theta-\vartheta_{j}\right),
\]
where $\gamma^{n_{Y}}$ is the order-$n_{Y}$ Dirichlet kernel
\[
\gamma^{2j+1}\left(\theta\right)=\sum_{k=-j}^{j}\mathrm{e}^{ik\theta}=\frac{\sin\left(\frac{2j+1}{2}\theta\right)}{\sin\frac{\theta}{2}}.
\]
In case $\boldsymbol{g}\left(\theta\right)=\theta+\omega$, where
$\omega$ is a constant rotation angle, we also have 
\[
\Omega_{ij}=\frac{1}{n_{Y}}\gamma^{n_{Y}}\left(\vartheta_{j}-\vartheta_{i}-\omega\right),
\]
which is a unitary matrix.

A higher dimensional rotation can also be calculated by using a rectangular
grid on $Y=\mathbb{T}^{d_{Y}}$ having $n_{1}\times n_{2}\times\cdots\times n_{d_{Y}}$
grid points, such that the elements of the function library are 
\[
\psi_{j_{1}\cdots j_{d_{Y}}}^{n_{1}\cdots n_{d_{Y}}}\left(\boldsymbol{\theta}\right)=\frac{1}{n_{1}\cdots n_{d_{Y}}}\gamma^{n_{1}}\left(\theta_{1}-\vartheta_{j}\right)\cdots\gamma^{n_{d_{Y}}}\left(\theta_{d_{Y}}-\vartheta_{j}\right).
\]
In case $\boldsymbol{g}\left(\theta\right)=\boldsymbol{\theta}+\boldsymbol{\omega}$
the $\boldsymbol{\Omega}$ matrix becomes
\[
\Omega_{\left(i_{1}\cdots i_{d_{Y}}\right)\left(j_{1}\cdots j_{d_{Y}}\right)}=\frac{1}{n_{1}\cdots n_{d_{Y}}}\gamma^{n_{1}}\left(\vartheta_{j_{1}}-\vartheta_{i_{1}}-\omega_{1}\right)\cdots\gamma^{n_{d_{Y}}}\left(\vartheta_{j_{d_{Y}}}-\vartheta_{i_{d_{Y}}}-\omega_{d_{Y}}\right).
\]

Parameter dependence on an interval $\theta\in Y=\left[a,b\right]$
can also be modelled, by choosing $\boldsymbol{\Psi}$ to be a polynomial
basis, such as 
\begin{equation}
\psi_{1}\left(\theta\right)=1,\psi_{2}\left(\theta\right)=\theta,\ldots,\psi_{n_{Y}}\left(\theta\right)=\theta^{n_{Y}-1},\label{eq:Parameter-Polynomial}
\end{equation}
or an appropriately scaled Chebyshev polynomial basis.

\subsection{\label{subsec:Bundles}Approximate stationary state and nearby linear
dynamics}

We can identify multiple invariant foliations from a single data set,
hence we need to make a choice which ones to recover. To make that
choice we fit an approximate linear map to the data. From the linear
map we calculate an approximate stationary state and the invariant
vector bundles about the stationary state, as defined in section \ref{subsec:Vector-Bundles}.
Once the vector bundles are calculated we select those that will form
the linear approximation of our invariant foliation as given by equations
\eqref{eq:Linear_ROM} and \eqref{eq:Linear_Encoder}.

Henceforth, we assume the following linear model

\begin{equation}
f_{i}\left(\boldsymbol{x},\boldsymbol{\alpha}\right)=A_{ijk}\alpha_{j}x_{k}+b_{ij}\alpha_{j}.\label{eq:LINID-model}
\end{equation}
The parameters $A_{ijk}$ and $b_{ij}$ of the model \eqref{eq:LINID-model}
are found by minimising the loss function
\begin{equation}
L_{lin}\left(\boldsymbol{A},\boldsymbol{b}\right)=\frac{1}{2}\sum_{t=0}^{N-1}\sum_{l=\ell_{t}}^{\ell_{t+1}-1}\left(A_{ijk}\alpha_{j,l}x_{k,l}+b_{ij}\alpha_{j,l}-x_{i,l+1}\right)^{2}.\label{eq:LINID-noscale}
\end{equation}

To find the stationary state, we need to solve the steady state equation
\eqref{eq:S_State_Invariance} for $\boldsymbol{s}$, which is written
as
\begin{equation}
s_{ij}\Omega_{jk}\alpha_{k}=A_{ijk}\alpha_{j}s_{kl}\alpha_{l}+b_{ij}\alpha_{j},\label{eq:LINID-S_State}
\end{equation}
when using the linear model \eqref{eq:LINID-model}. Rearranging equation
\eqref{eq:LINID-S_State} yields the linear system
\[
\left(\delta_{ik}\Omega_{jl}-A_{i\underline{l}k}\delta_{j\underline{l}}\right)s_{kj}=b_{il}
\]
that is solved for the matrix $s_{kj}$. Now we define the matrix
\[
B_{\left(il\right)\left(kj\right)}=\delta_{ik}\Omega_{jl}-A_{i\underline{l}k}\delta_{j\underline{l}}
\]
and by inverting this matrix we find that the solution to \eqref{eq:LINID-S_State}
is
\[
s_{kj}=\left[B^{-1}\right]_{kjil}b_{il}.
\]
The next task is to find the invariant vector bundles using equation
\eqref{eq:ED-left-bundle}. Numerically we can make the matrices $\Lambda_{j}$
constant, and therefore we solve the eigenvalue-eigenvector problem
\begin{equation}
\lambda_{p}U_{p}\left(\boldsymbol{\theta}\right)=U_{p}\left(\boldsymbol{g}\left(\boldsymbol{\theta}\right)\right)\boldsymbol{A}\left(\boldsymbol{\theta}\right),\label{eq:LINID-eig}
\end{equation}
where $\lambda_{p}$ are the diagonal elements of the diagonal matrices
$\boldsymbol{\Lambda}_{j}$. The eigenvalue problem \eqref{eq:LINID-eig}
in discrete form is
\[
\lambda_{p}U_{pit}\alpha_{i}=U_{pij}\Omega_{ik}\alpha_{k}A_{jlt}\alpha_{l},
\]
which after rearrangement and eliminating $\alpha$ becomes
\begin{equation}
U_{pij}\left(\lambda_{p}\delta_{jt}\delta_{il}-\Omega_{i\underline{l}}A_{j\underline{l}t}\right)=0.\label{eq:LINID-eigen-problem}
\end{equation}

There are $n_{Y}$-times as many eigenvalues and eigenvectors than
the number of one-dimensional vector bundles. This means that for
each vector bundle there should be $n_{Y}$ eigenvectors. We create
an adjacency matrix between all pairs of eigenvectors
\[
\mathit{Adj}_{pq}=\sum_{i=1,j=1}^{n_{Y},d_{X}}\left|\nu_{qi}U_{pij}-\nu_{pi}U_{qij}\right|^{2},
\]
where 
\[
\nu_{pi}=\sqrt{U_{\underline{p}\underline{i}j}\overline{U}_{\underline{p}\underline{i}j}}
\]
and find $d_{X}$ clusters of $n_{Y}$ eigenvectors that are closest
to each other. Within each cluster the eigenvector with smallest total
variation is chosen, where total variation is calculated by the formula
\[
\mathit{TV}_{p}=\sum_{i=2,j=1}^{n_{Y},d_{X}}\left|U_{pij}-U_{p\left(i-1\right)j}\right|.
\]
This yields a set of $d_{X}$ indices $p_{1},p_{2},\ldots,p_{d_{X}}$.
The indices are sorted in such a way that up until index $p_{c}$
we have complex conjugate pairs
\[
\lambda_{p_{i}}=\overline{\lambda}_{p_{i+1}},\;i<c
\]
 and 
\[
\left|\lambda_{p_{1}}\right|\ge\left|\lambda_{p_{2}}\right|\ge\cdots\ge\left|\lambda_{p_{c}}\right|,\quad\left|\lambda_{p_{c+1}}\right|\ge\left|\lambda_{p_{c+2}}\right|\ge\cdots\ge\left|\lambda_{p_{n_{X}}}\right|.
\]
Then we have overall $m_{\Sigma}=d_{X}-c/2$ real vector bundles.
To make sure that the data is well-scaled when it is encoded with
the linear encoders $U_{j}$, we also need to scale the linear maps
$\boldsymbol{U}_{j}$ with respect to the data. The data components
then become
\begin{align*}
\hat{x}_{2j-1,l} & =\Re U_{p_{2j-1}ik}\alpha_{i,l}x_{k,l} & \text{for }1\le j\le\frac{c}{2},\\
\hat{x}_{2j,l} & =\Im U_{p_{2j-1}ik}\alpha_{i,l}x_{k,l} & \text{for }1\le j\le\frac{c}{2},\\
\hat{x}_{j,l} & =U_{p_{j-c/2}ik}\alpha_{i,l}x_{k,l} & \text{for }\frac{c}{2}<j\le d_{X}-c/2.
\end{align*}
With this transformed data we calculate the normalising factors
\[
\nu_{i}=\sqrt{d_{X}}\max_{j}\left|\hat{x}_{i,j}\right|.
\]
The linear parts of the model and the encoder, as appear in \eqref{eq:Linear_ROM}
and \eqref{eq:Linear_Encoder}, respectively become
\begin{align*}
\breve{\boldsymbol{U}}_{j}\left(\boldsymbol{\alpha}\right)\boldsymbol{x} & =\begin{pmatrix}\begin{array}{l}
\nu_{2j-1}^{-1}\Re U_{p_{2j-1}ik}\alpha_{i}x_{k}\\
\nu_{2j}^{-1}\Im U_{p_{2j-1}ik}\alpha_{i}x_{k}
\end{array}\end{pmatrix}, & \boldsymbol{\Lambda}_{j} & =\begin{pmatrix}\Re\lambda_{p_{2j-1}} & -\frac{\nu_{2j-1}}{\nu_{2j}}\Im\lambda_{p_{2j-1}}\\
\frac{\nu_{2j}}{\nu_{2j-1}}\Im\lambda_{p_{2j-1}} & \Re\lambda_{p_{2j-1}}
\end{pmatrix}, & 1\le j\le\frac{c}{2},\\
\breve{\boldsymbol{U}}_{j}\left(\boldsymbol{\alpha}\right)\boldsymbol{x} & =\nu_{j-c/2}^{-1}U_{p_{j-c/2}ik}\alpha_{i}x_{k}, & \boldsymbol{\Lambda}_{j} & =\lambda_{p_{j-c/2}}, & \frac{c}{2}<j\le m_{\Sigma}.
\end{align*}
The overall linear transformation becomes
\[
\breve{\boldsymbol{U}}\left(\boldsymbol{\alpha}\right)=\begin{pmatrix}\breve{\boldsymbol{U}}_{1}\left(\boldsymbol{\alpha}\right)\\
\vdots\\
\breve{\boldsymbol{U}}_{m_{\Sigma}}\left(\boldsymbol{\alpha}\right)
\end{pmatrix}
\]
and its inverse is
\[
\breve{\boldsymbol{W}}\left(\boldsymbol{\alpha}\right)=\breve{\boldsymbol{U}}\left(\boldsymbol{\alpha}\right)^{-1}.
\]
For a given index set we create two sets of transformed data. The
first has the same dimension as the original data set
\[
\breve{\boldsymbol{x}}_{j}=\breve{\boldsymbol{U}}\left(\boldsymbol{\alpha}_{j}\right)\boldsymbol{x}_{j},\quad j=1,\ldots,\ell_{N}.
\]
The second set contains a subset of the coordinates of $\breve{\boldsymbol{x}}$,
which do no belong to the index set $\mathcal{I}$, that is 
\begin{equation}
\boldsymbol{x}_{j}^{\perp}=\breve{\boldsymbol{U}}^{\perp}\left(\boldsymbol{\alpha}\right)\boldsymbol{x}_{j}=\begin{pmatrix}\breve{\boldsymbol{U}}_{j_{1}}\left(\boldsymbol{\alpha}\right)\\
\vdots\\
\breve{\boldsymbol{U}}_{j_{\bar{p}}}\left(\boldsymbol{\alpha}\right)
\end{pmatrix}\boldsymbol{x}_{j},\label{eq:LINID-Uperp}
\end{equation}
where $\left\{ j_{1},\ldots,j_{\bar{p}}\right\} =\left\{ 1,\ldots,m_{\Sigma}\right\} \setminus\mathcal{I}$.
Note that $\breve{\boldsymbol{U}}^{\perp}\left(\boldsymbol{\alpha}\right)$
is implicitly defined by equation \eqref{eq:LINID-Uperp}.

\subsection{\label{subsec:Encoder-representation}Encoder representation}

\begin{figure}
\begin{centering}
\includegraphics[width=0.45\textwidth]{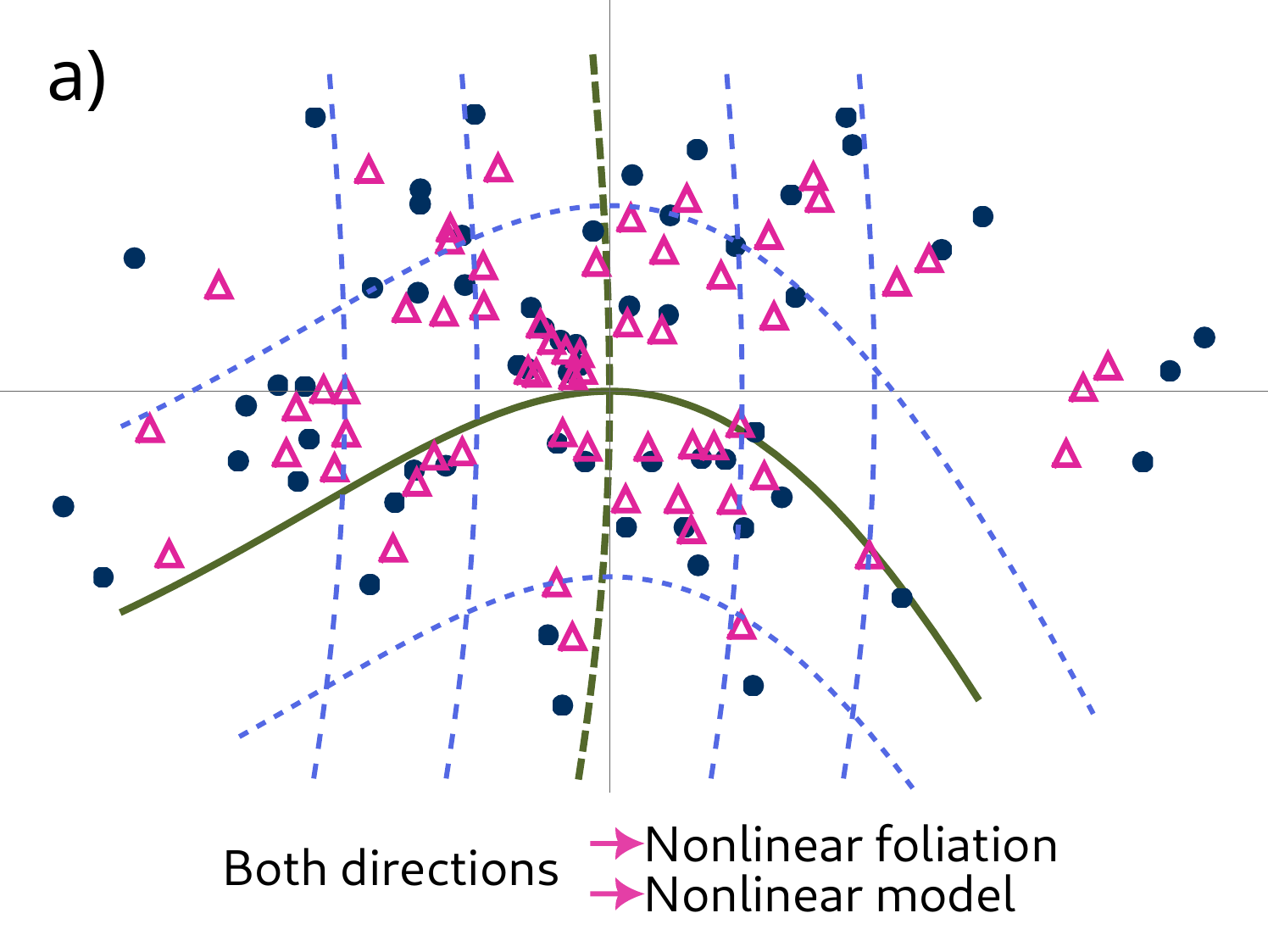}\includegraphics[width=0.45\textwidth]{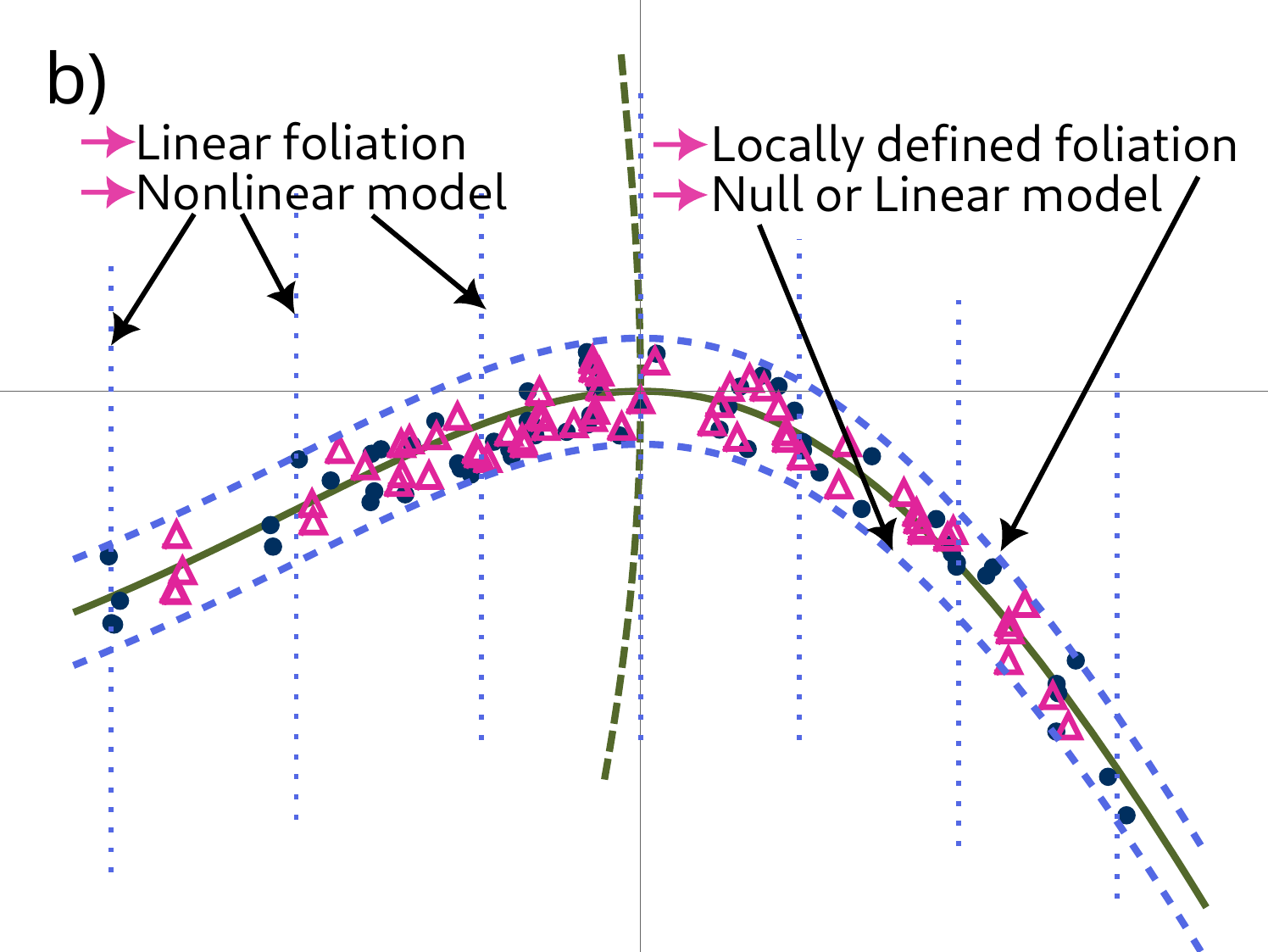}
\par\end{centering}
\caption{\label{fig:Data-Distribution}Two typical data distributions. a) The
data is evenly distributed in phase space and therefore it is appropriate
to use a nonlinear encoder and a nonlinear model for both possible
foliations. b) The data is distributed in the neighbourhood of an
invariant manifold. The invariant foliation that resolves the dynamics
on the invariant manifold can be represented by a linear foliation
and a nonlinear model. The invariant foliation representing the transversal
dynamics can be represented by a null or linear model and the encoder
can be locally defined.}
\end{figure}
The encoder is represented by a polynomial with low-rank tensor coefficients
\cite{TensorApproxSurvey}. The first representation does not assume
reducibility, hence generally applicable. The representation is

\begin{equation}
u_{i}\left(\breve{\boldsymbol{x}},\boldsymbol{\alpha}\right)=u_{ij}^{0}\alpha_{j}+u_{ijk}^{1}\alpha_{j}\breve{x}_{k}+u_{ijk_{1}}^{2}\alpha_{j}x_{k_{1}}^{\perp}\breve{x}_{k_{2}}+\cdots+u_{ijk_{1}\cdots k_{\mathit{EO}}}^{\mathit{EO}}\alpha_{j}x_{k_{1}}^{\perp}\breve{x}_{k_{2}}\cdots\breve{x}_{k_{\mathit{EO}}}\label{eq:ENC-generic}
\end{equation}
with the assumption that 
\begin{equation}
u_{i_{1}j\underline{k}}^{1}u_{i_{2}j\underline{k}}^{1}=\delta_{i_{1}i_{2}},\label{eq:ENC-norm-generic}
\end{equation}
where the summation does not apply over index $k$, meaning that the
$u_{ijk}^{1}$ is an orthogonal matrix for each fixed index $k$.
The superscript $^{\mathit{EO}}$ represents the encoder order. The
constraint \eqref{eq:ENC-norm-generic} makes sure that the encoder
does not collapse. Moreover, the use of $\boldsymbol{x}^{\perp}$
in each nonlinear term makes sure that the encoder is linear for the
vector bundle where $\boldsymbol{x}^{\perp}=0$. The vector bundle
for which the encoder is linear can also be written as 
\[
\mathcal{B}^{\parallel}=\left\{ \breve{\boldsymbol{U}}\left(\boldsymbol{\alpha}\right)\mathrm{null}\,\breve{\boldsymbol{U}}^{\perp}\left(\boldsymbol{\alpha}\right):\boldsymbol{\alpha}\in\boldsymbol{\Psi}\left(Y\right)\right\} .
\]
The linearity constraint is analogous to the graph-style parametrisation
of invariant manifolds \cite{CabreLlave2003}. One disadvantage of
representation \eqref{eq:ENC-generic} is that it has a large number
of parameters. Figure \ref{fig:Data-Distribution} illustrates two
possible data distributions. Representation \eqref{eq:ENC-generic}
is suitable for the type of data in \ref{fig:Data-Distribution}(a).
If the data is clustered about a manifold, the parameters of \eqref{eq:ENC-generic}
become ill-defined. If the sought after foliation has leaves transversal
to the manifold, the encoder can be made linear. If the desired foliation
has leaves aligned with the manifold, a locally defined encoder should
be used. This scenario is illustrated in \ref{fig:Data-Distribution}(b).

A locally defined encoder has less parameters than \eqref{eq:ENC-generic}.
The linear part of the encoder remains pointwise orthogonal, constrained
by equation \eqref{eq:ENC-norm-generic}, but the nonlinear part depends
only on $\boldsymbol{x}^{\perp}$ which has the same dimensionality
as a leaf within the foliation. The locally defined encoder can be
written as
\begin{equation}
u_{i}\left(\breve{\boldsymbol{x}},\boldsymbol{\alpha}\right)=u_{ij}^{0}\alpha_{j}+u_{ijk}^{1}\alpha_{j}\breve{x}_{k}+u_{ijk_{1}}^{2}\alpha_{j}x_{k_{1}}^{\perp}x_{k_{2}}^{\perp}+\cdots+u_{ijk_{1}\cdots k_{\mathit{EO}}}^{\mathit{EO}}\alpha_{j}x_{k_{1}}^{\perp}\cdots x_{k_{\mathit{EO}}}^{\perp}.\label{eq:ENC-local}
\end{equation}
The level surfaces of \eqref{eq:ENC-local} are shifted versions of
each other and therefore \eqref{eq:ENC-local} is useful if the data
is clustered about a single level surface (typically, the zero-level
surface) of \eqref{eq:ENC-local}.

If the system is reducible, the conjugate map is autonomous and encodes
less information. In this case the encoder must be more general than
\eqref{eq:ENC-generic} and must allow for a variation in the norm
of the matrix $u_{ijk}^{1}\alpha_{j}$, while making sure that the
tensor $u_{ijk}^{1}$ does not vanish. The form for this encoder representation
is
\begin{equation}
u_{i}\left(\breve{\boldsymbol{x}},\boldsymbol{\alpha}\right)=u_{ij}^{0}\alpha_{j}+u_{ijk}^{1}\alpha_{j}\breve{x}_{k}+u_{ijk_{1}}^{2}\alpha_{j}\breve{x}_{k_{1}}\breve{x}_{k_{2}}+\cdots+u_{ijk_{1}\cdots k_{\mathit{EO}}}^{\mathit{EO}}\alpha_{j}\breve{x}_{k_{1}}\cdots\breve{x}_{k_{\mathit{EO}}}.\label{eq:ENC-reducible}
\end{equation}
To prevent collapse, we constrain the linear part of \eqref{eq:ENC-reducible}
in a average sense over the full vector bundle, that is, we stipulate
that matrix $u_{i\left(jk\right)}^{1}$ is orthogonal. In formulae
the constraint is 
\begin{equation}
u_{i_{1}jk}^{1}u_{i_{2}jk}^{1}=\delta_{i_{1}i_{2}},\label{eq:ENC-norm-local}
\end{equation}
where the summation also applies over index $k$.

In practice it is possible to combine the locally defined encoder
\eqref{eq:ENC-local}, with the assumption of reducibility and apply
the constraint \eqref{eq:ENC-norm-local} instead of \eqref{eq:ENC-norm-generic}.
In case of a locally defined foliation, the conjugate map can be linear
or even null, that is $\boldsymbol{r}\left(\boldsymbol{z},\boldsymbol{\theta}\right)=\boldsymbol{0}$.
This is because there might be very little or no information in the
dynamics transversal to the invariant manifold.

\subsection{Summary of the procedure}

In order to obtain reduced order models the following steps can be
taken.
\begin{enumerate}
\item If necessary, apply state-space reconstruction as described in section
\ref{subsec:State-space-reconstruction}.
\item Find an approximate linear model and the steady state using the procedure
in section \ref{subsec:Bundles}. This is equivalent to dynamic mode
decomposition, which removes redundancy from the data and brings it
to a coordinate system, where the approximate linear model is block
diagonal.
\item Choose a set of invariant vector bundles of the linear model and fit
a foliation to them as in section \ref{subsec:Fitting-a-foliation}.
The conjugate map can be a polynomial, and a choices of encoder representations
are described in section \ref{subsec:Encoder-representation}.
\item Extract invariant manifolds and ascertain the accuracy of the invariant
foliations using testing data. The invariant manifold can be found
using the method in section \ref{subsec:Foliations-to-Manifolds}.
In the special case of a two-dimensional invariant manifold, backbone
curves can be recovered by the method presented in section \ref{subsec:Normal-form}.
Note that the backbone curves fully specify the reduced order model,
because the polar models \eqref{eq:Polar-Model-Map} and \eqref{eq:Polar-Model-ODE}
can be uniquely reproduced from them using the relations \eqref{eq:INST-damp},
\eqref{eq:INST-freq} or \eqref{eq:INST-damp-ode}, \eqref{eq:INST-freq-ode}
for differential equations.
\end{enumerate}

\section{Examples}

Here we demonstrate how the theory and its numerical implementation
works through a variety of examples. The examples include autonomous,
parameter dependent and periodically forced systems. Some of the examples
use experimental data, others use synthetic data. The examples also
illustrate how various kinds of data can be used, ranging from a single
trajectory near an invariant manifold to high-dimensional systems
with data far from any invariant manifold.

\subsection{Titanium beam}

The first example is a single trajectory near an invariant manifold.
The experimental setup consists of a titanium beam suspended by two
springs. The experiment was initially designed to measure gravitational
interaction between two beams \cite{Brack2022}. The data we use was
published alongside \cite{Bettini2025}. 
\begin{figure}
\begin{centering}
\includegraphics[width=0.25\textwidth]{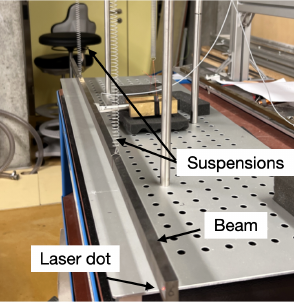}\includegraphics[width=0.49\textwidth]{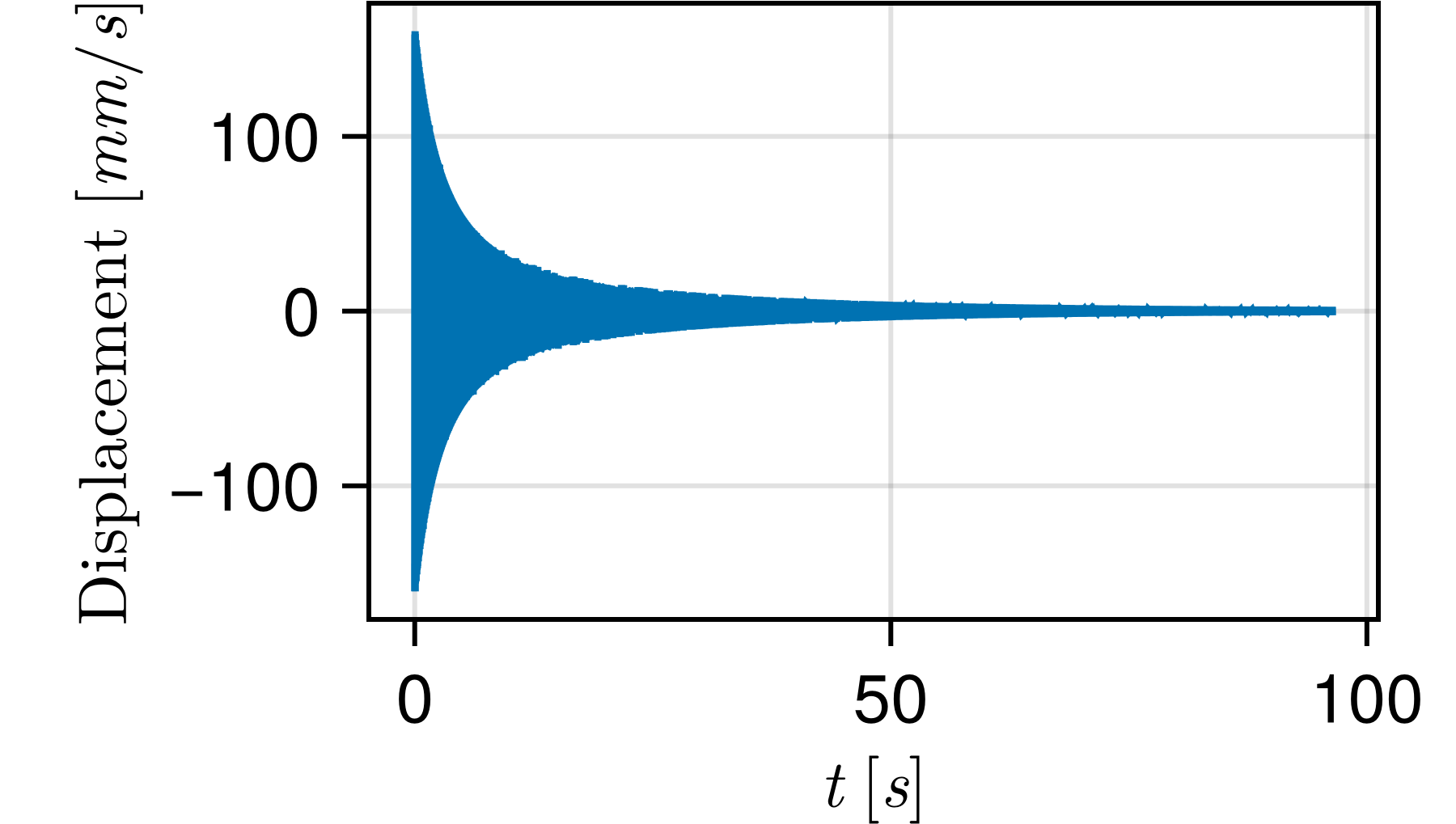}
\par\end{centering}
\caption{\label{fig:GRAV-rig}Left: experimental setup, reproduced from \cite{Bettini2025}.
Right: illustration of the available data.}
\end{figure}

The data is a single trajectory with a scalar output measuring the
velocity of a single point on a titanium beam. The initial condition
of the beam was set by forcing the beam at its first resonant frequency,
then stopping the forcing and recording the decaying vibration. This
set-up procedure is required by the SSMLearn method \cite{Bettini2025}
to ensure that the data is close to an invariant manifold. The present
method does not require setting a precise initial condition, impact
hammer tests are sufficient, as shown by other examples below. The
damping of the system is extremely low, hence the trajectory is long
and represents a significant amount of data. Test data is not available,
so it is impossible to check if the model generalises to unseen data.

The signal is sampled with period $\Delta t=0.13653$ ms. We use $d=64$
delay length to reconstruct the phase space, which is then used to
identify an approximate linear model. The first 6 dynamic modes are
used for constructing two invariant foliations. The rest of the dynamic
modes are within the noise floor and contain no useful information.
The first invariant foliation aligns with dynamic modes 1,2, which
has a linear encoder and an order-7 conjugate map. The second invariant
foliation takes dynamic modes 3-6 and uses an order-7 locally defined
encoder and a linear conjugate map. This is conceptually similar to
the scenario illustrated in figure \ref{fig:Data-Distribution}(b).

\begin{figure}
\begin{centering}
\includegraphics[width=0.99\textwidth]{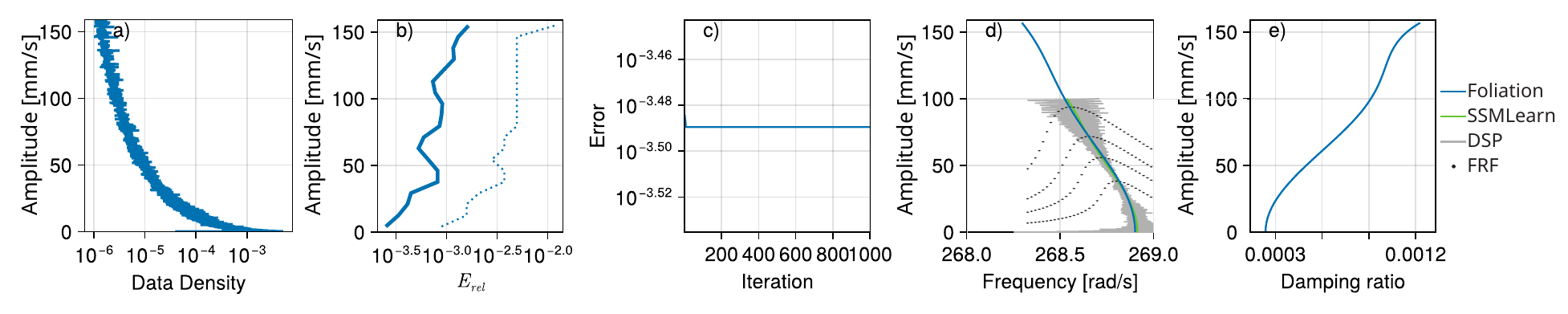}
\par\end{centering}
\caption{\label{fig:GRAV-result}Results of fitting an invariant foliation
to the data. a) Data density as a function of vibration amplitude.
b) Solid line represents the mean relative error $E_{\mathit{rel}}$
as given by equation \eqref{eq:Relative_Error}, the dotted line is
maximum relative error at a given amplitude. c) Training error as
a function of optimisation steps. The method converges within the
first few steps. d) Instantaneous frequency of the first vibration
mode as calculated from the invariant foliation (blue line) and compared
to the result from \cite{Bettini2025}. The line labelled DSP is calculated
using the method from \cite{JinBrake2020FreqDamp} by the authors
of \cite{Bettini2025}. The dots represent experimental forced response.
e) Instantaneous damping as defined by equation \eqref{eq:INST-damp-wrong}.}

\end{figure}

The result can be seen in figure \ref{fig:GRAV-result}. The optimisation
converged quickly to the critical point. The instantaneous frequency
matches the analysis in \cite{Bettini2025}. Moreover, our method
was able to identify the backbone curves for the whole of the trajectory,
which is 50\% more than in \cite{Bettini2025}. This is because SSMLearn
requires discarding data that is not close to an invariant manifold,
thereby losing significant information, while our method uses all
data and the encoder accurately picks out the required vibration mode.

\begin{figure}
\begin{centering}
\includegraphics[width=0.35\textwidth]{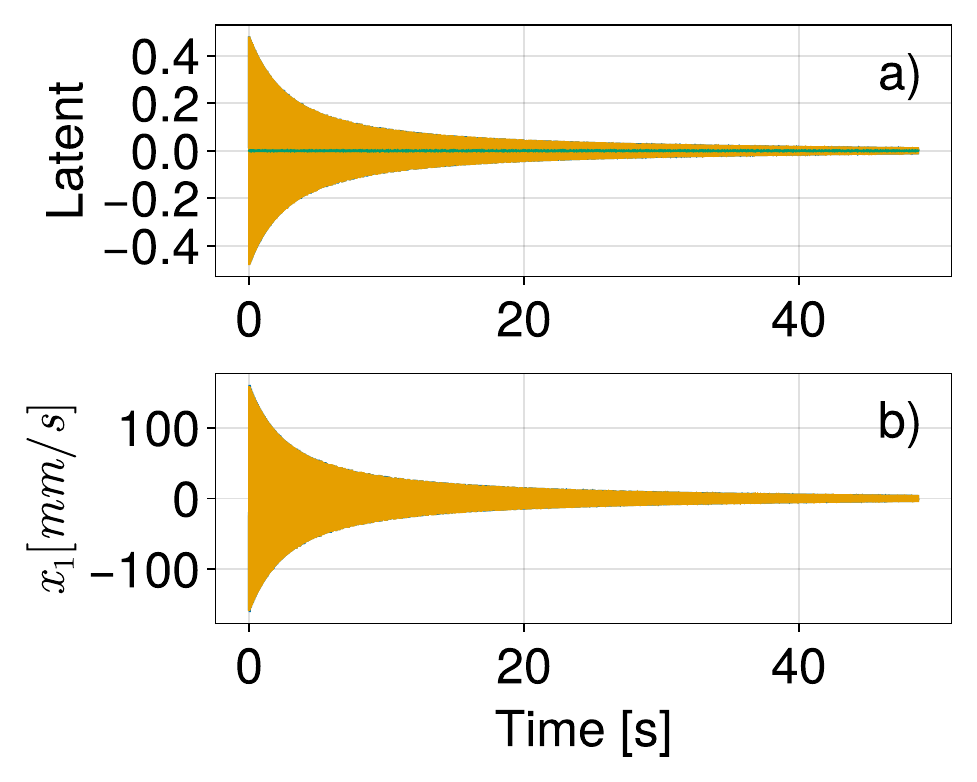}
\par\end{centering}
\caption{\label{fig:GRAV-latent}Error in latent space. The blue line (covered)
is the encoded trajectory, the mustard line is the model prediction
and the green line is the difference of the two. The diagram demonstrates
that the model and data are fully synchronised in phase and follow
the same envelope without noticeable error.}

\end{figure}

The accuracy of the method is measured in latent space, which means
that iterations of the conjugate map with an appropriate initial condition
must match the encoded data. In case of this example the error is
negligible throughout the whole of the trajectory as illustrated in
figure \ref{fig:GRAV-latent}(a). Due to the experimental setup there
is unnoticeable signal in the second invariant foliation, hence the
model prediction in physical space also matches the data perfectly,
as displayed in figure \ref{fig:GRAV-latent}(b).

\subsection{Shaw-Pierre model}

\begin{figure}
\begin{centering}
\includegraphics{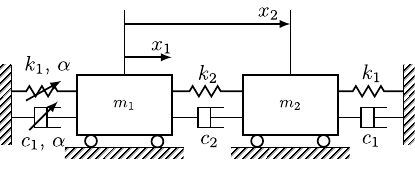}
\par\end{centering}
\caption{\label{fig:Shaw-Pierre}Schematic of Shaw-Pierre oscillator as presented
in \cite{Bettini2025}.}

\end{figure}
The Shaw-Pierre model \cite{ShawPierre} is a classical example demonstrating
that invariant manifolds are suitable nonlinear normal modes. This
example was used in many subsequent papers, including \cite{Bettini2025}
as a benchmark. The model features minimal nonlinearity, in the form
of a nonlinear spring and a nonlinear damper. Here, we also add forcing
to the system to fully demonstrate our method. The schematic of the
system can be seen in figure \ref{fig:Shaw-Pierre}. The governing
equations of the system are

\begin{equation}
\dot{\boldsymbol{x}}=\begin{pmatrix}0 & 0 & 1 & 0\\
0 & 0 & 0 & 1\\
-k_{1} & k_{2} & -c_{1} & c_{2}\\
0 & -\left(k_{1}+k_{2}\right)+k_{2} & 0 & -\left(c_{1}+c_{2}\right)+c_{2}
\end{pmatrix}\boldsymbol{x}+\begin{pmatrix}0\\
0\\
-\alpha x_{1}^{3}+\alpha c_{1}x_{3}^{3}+\beta\cos(\omega(t+0.1))\\
\beta\cos t
\end{pmatrix}\label{eq:SP-model}
\end{equation}
and parameters values are
\[
k_{1}=1,\;k_{2}=3.325,\;c_{1}=0.05,\;c_{2}=0.01,\;\alpha=0.5.
\]
We have created 12 training trajectories with sampling period $\Delta t=0.2780$
s. We have also created a testing trajectory to assess the accuracy
of the invariant foliation. The results are discussed for the invariant
foliation of the first natural frequency only, because the invariant
foliation for the second natural frequency is linear. Due to the duality
of foliations and manifolds, the invariant manifold of the first natural
frequency is a linear subspace with nonlinear dynamics on it, which
makes it trivial to calculate using SSMLearn \cite{Bettini2025}.

The first invariant foliation is assumed to be a cubic polynomial
with a cubic encoder, the second invariant foliation is assumed to
be linear both for the conjugate map and the encoder. For the autonomous
and the forced system, the encoder of the first natural frequency
are of type \eqref{eq:ENC-reducible}, while for the linear encoder
of the second natural frequency the encoder is of type \eqref{eq:ENC-generic}.
The periodic forcing is represented at $n_{Y}=19$ collocation points
as defined by \eqref{eq:circle-colloc}. For the parametric invariant
foliation all encoders are of type \eqref{eq:ENC-generic}, because
parameter dependence is not reducible.

\begin{figure}
\begin{centering}
\includegraphics[width=0.99\textwidth]{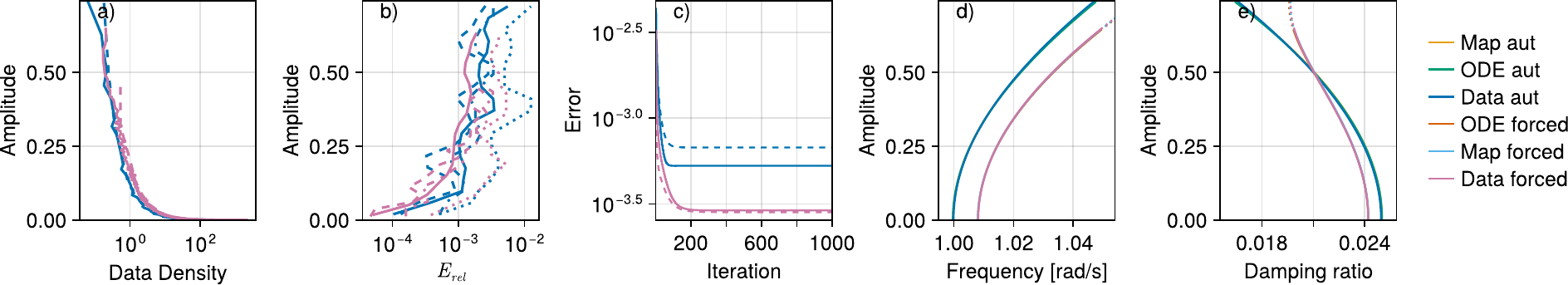}
\par\end{centering}
\caption{\label{fig:SP-autonomous}First natural frequency of the model in
figure \ref{fig:Shaw-Pierre}. The blue and purple lines represent
the autonomous and the forced system, respectively a) Data density:
continuous lines denote the training data dashed lines denote testing
data. b) Relative error \eqref{eq:Relative_Error}: continuous lines
are mean relative error of training data, dotted lines are maximum
relative error of training data, dashed lines are mean relative error
of testing data, dash-dotted lines are maximum relative error of testing
data. c) Continuous lines are the training errors, dashed lines are
the testing errors as functions of optimisation steps. d) Frequency
backbone curve as given by equation \eqref{eq:INST-freq}, dots mean
extrapolation, without training data. e) Damping ratio curve as given
by equation \eqref{eq:INST-damp}, dots mean extrapolation. The backbone
curves in d) and e) were also calculated from the differential equation
\eqref{eq:SP-model} and its time-$\Delta t$ solution map using remarks
\ref{rem:ODE-backbones} and \ref{rem:MAP-backbones}, respectively.}
\end{figure}

The analysis of the identified invariant foliations can be seen in
figure \ref{fig:SP-autonomous} both for the autonomous ($\beta=0$)
and the forced case with $\beta=0.04$, $\omega=1.1892$ rad/s. Forcing
has slightly increased the frequency of the free decaying signal about
the steady state. It can also be seen that the invariant foliation
extrapolates beyond the vibration amplitudes present in the training
data, because the extrapolated backbone curves match the directly
calculated curves from the model equation \eqref{eq:SP-model}.

\begin{figure}
\begin{centering}
\includegraphics[width=0.49\textwidth]{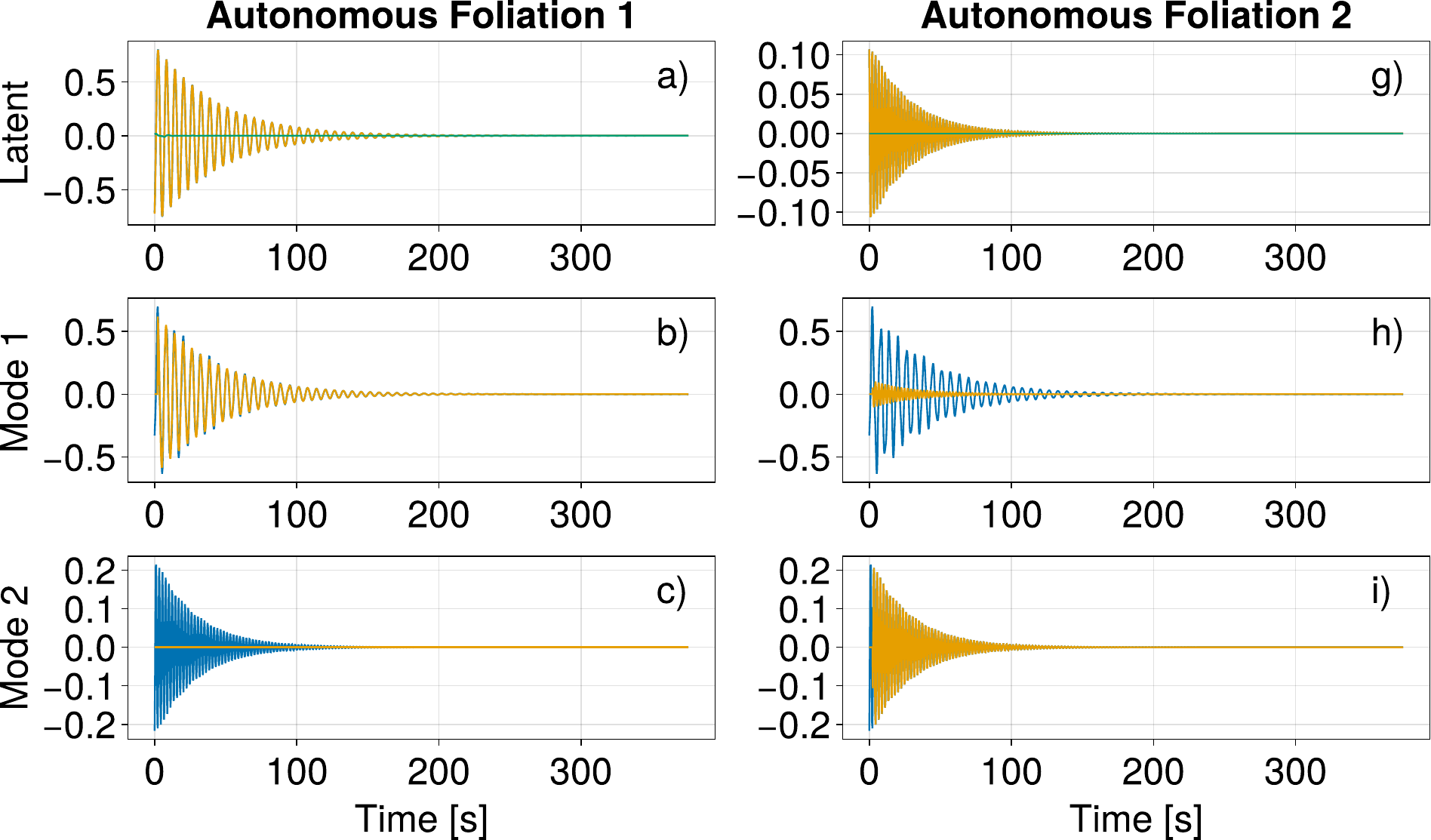}\includegraphics[width=0.49\textwidth]{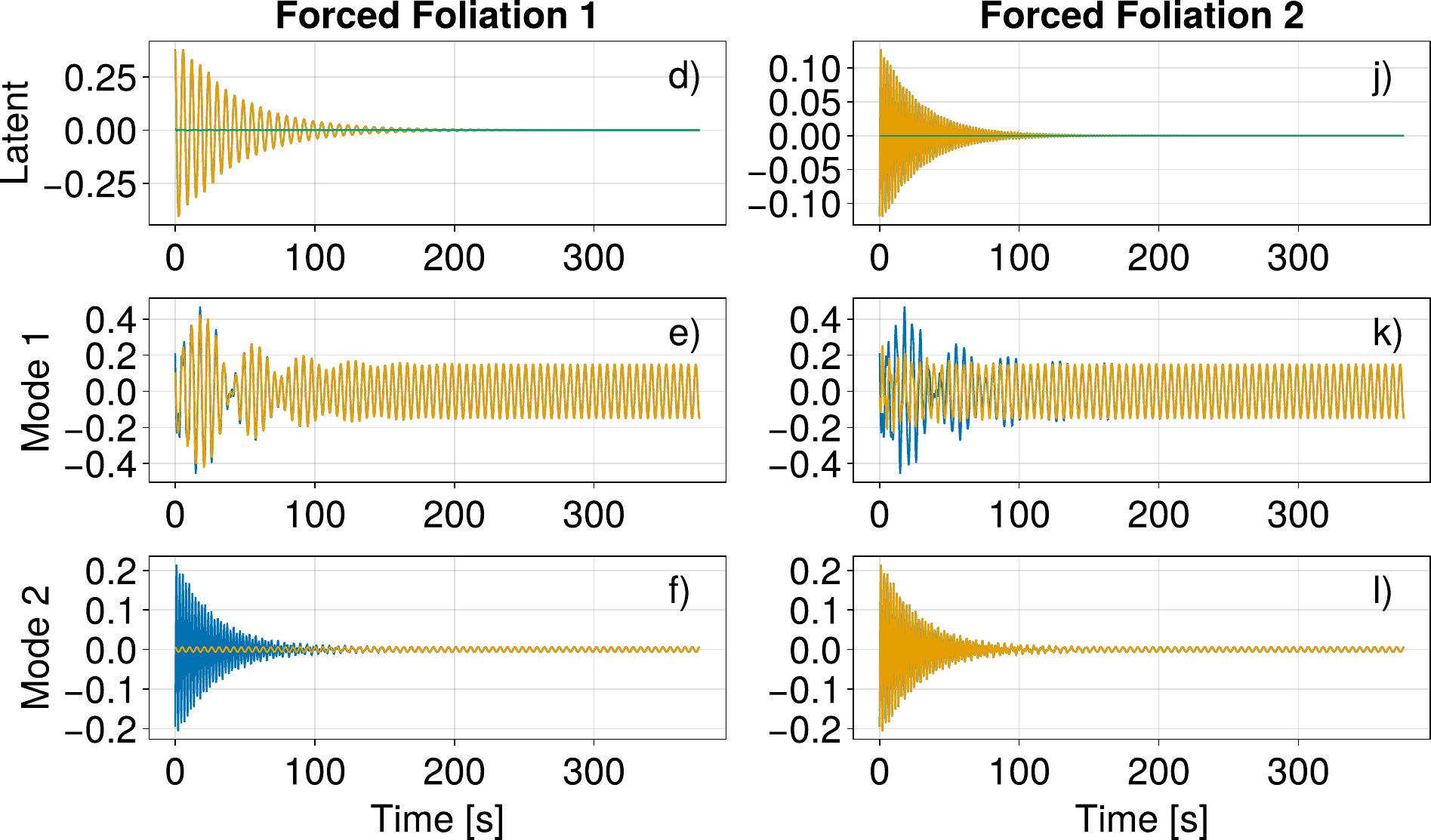}
\par\end{centering}
\caption{\label{fig:SP-latent}Visualising the accuracy of the invariant foliations
on the testing trajectory for the model in figure \ref{fig:Shaw-Pierre}.
The first two columns (left to right) and the last two columns are
for autonomous and forced systems, respectively. The top row displays
one of (out of the two) the latent variables, subsequent rows display
the displacements $x_{1}$ and $x_{2}$. The blue lines are the testing
trajectory, the mustard lines are the model prediction and the green
lines represent the error.}
\end{figure}

The accuracy of the calculation can be visualised by comparing trajectories
in latent space. Figures \ref{fig:SP-latent}(a,d,g,j) show that the
identified models predict the encoded trajectories with high accuracy.
The model prediction can also be mapped back to physical coordinates.
In figures \ref{fig:SP-latent}(b,c,e,f,h,i,k,l) we compare the displacements
of the two masses ($x_{1},x_{2}$) to the reconstructed model trajectories.
Due to linearity of the second invariant foliation and the way the
mechanical system is modelled, the displacement between the two masses
is independent of the first vibration mode. 

To identify parameter dependent foliations, parameter $\alpha$ was
varied in the range of $\left[0.4,0.6\right]$ on a Chebyshev grid
of 16 points, resulting in 16 training trajectories with different
parameters. The parameter dependence of all quantities was modelled
using quartic polynomial interpolation as given by the basis in equation
\eqref{eq:Parameter-Polynomial}. A single testing trajectory is produced
at $\alpha=0.48988$ , which is not in the training set. All quantities
were then calculated for the parameter value of the testing trajectory
and visualised in figure \ref{eq:Parameter-Polynomial}. The accuracy
is similar to all previous calculations.

\begin{figure}
\begin{centering}
\includegraphics[width=0.99\textwidth]{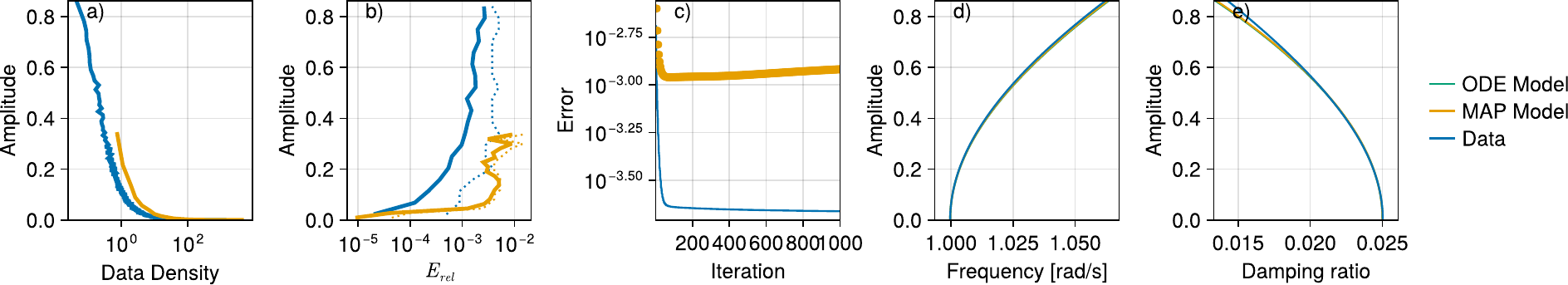}
\par\end{centering}
\caption{\label{fig:SP-parametric}The same type of result as in figure \ref{fig:SP-autonomous},
but the system is trained with different $\alpha$ parameter values
from the interval $\left[0.4,0.6\right]$. The result is only plotted
for the parameter of the testing trajectory with parameter $\alpha=0.48988$,
which is not represented in the training set.}
\end{figure}

\subsection{Car-following model}

This example shows that normal hyperbolicity can break down along
an invariant manifold, which is then picked up by direct numerical
calculations of the invariant manifold. Surprisingly, this breakdown
has little effect when the invariant manifold is extracted from the
data-driven invariant foliation. 

\begin{figure}
\begin{centering}
\includegraphics[width=0.35\textwidth]{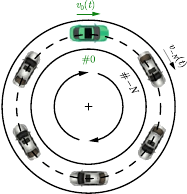}
\par\end{centering}
\caption{\label{fig:Carring}Five cars follow each other on a circular track.
The corresponding mathematical model is equation \eqref{eq:CarFollow-Model}.}

\end{figure}
The schematic of the system can be seen in figure \ref{fig:Carring},
where $n=5$ cars follow each other. A detailed description of the
model can be found in \cite{OroszWilsonEtAl2009}. The cars on the
track cannot overtake each other and their velocity $v_{k}$ is dictated
by the distance from the car just in front of them $h_{k}$ (called
headway). In this particular model drivers react instantly to the
changes in headway. Each car has a maximum velocity $V_{k}$. The
differential equations describing the system are

\begin{align}
\dot{v}_{k} & =\alpha\left(\frac{V_{k}\left(h_{k}-1\right)^{2}}{1+\left(h_{k}-1\right)^{2}}-v_{k}\right), &  & k=1\ldots n,\nonumber \\
\dot{h}_{k} & =v_{k-1}-v_{k}, &  & k=2\ldots n,\label{eq:CarFollow-Model}\\
h_{1} & =L-\sum_{j=2}^{n}h_{j}.\nonumber 
\end{align}
We assume that the length of the circular track is $L=2n$, the maximum
velocity of each vehicle is $V_{k}=1$, except for $V_{n}=1+A\cos\omega t$,
which is time-dependent. We set the forcing frequency to $\omega=0.63246$
rad/s and use either $A=0$ or $A=0.1$ as the forcing amplitude.
For this given set of parameters and without forcing, the steady state
is such that all headways are equal $h_{k}^{\star}=L/n$ and all cars
travel at speed 
\[
v_{k}^{\star}=\frac{\left(L/n-1\right)^{2}}{1+\left(L/n-1\right)^{2}}.
\]
The value $\alpha=0.75$, that controls the acceleration of a car,
is used so that the system is stable but near the stability boundary.
The eigenvalues of the Jacobian about this equilibrium are
\begin{align*}
\lambda_{12} & =-0.0163\pm0.4971i,\\
\lambda_{34} & =-0.2276\pm0.7480i,\\
\lambda_{56} & =-0.5223\pm0.7480i,\\
\lambda_{78} & =-0.7337\pm0.4971i,\\
\lambda_{9} & =-0.75.
\end{align*}
The spectral quotient for first two eigenvalues is $\beth_{1-2}=1$
while the spectral quotient for the remaining eigenvalues is $\beth_{3-9}=45.993$.
To discretise the system we used $n_{Y}=17$ collocation points for
the forcing. The forced system about the periodic orbit has Floquet
exponents
\begin{align*}
\lambda_{12} & =-0.0209\pm0.5040i,\\
\lambda_{34} & =-0.2281\pm0.7409i,\\
\lambda_{56} & =-0.5203\pm0.7409i,\\
\lambda_{78} & =-0.7306\pm0.4919i,\\
\lambda_{9} & =-0.75.
\end{align*}
The spectral quotients are $\beth_{1-2}=1$ and $\beth_{3-9}=33.7454$.
The high value for the spectral quotient means that the dynamics for
that foliation will decay to the steady state very quickly.

\begin{figure}
\begin{centering}
\includegraphics[width=0.99\textwidth]{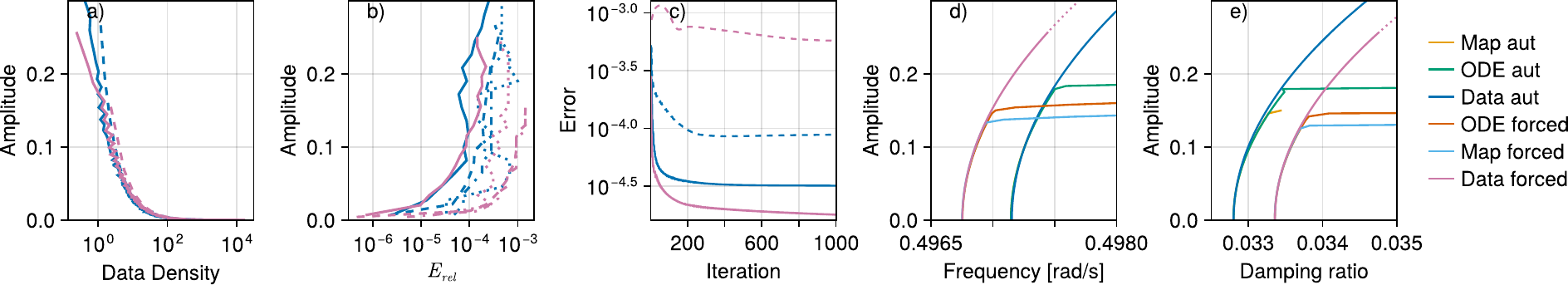}
\par\end{centering}
\caption{\label{fig:CarFollow-Auto}The reduced order model and its diagnostics
for the first mode of the car following model \eqref{eq:CarFollow-Model}.
The interpretation of the curves is the same as in figure \ref{fig:SP-autonomous}.}
\end{figure}

In all cases two invariant foliations are calculated, one corresponding
to eigenvalues $\lambda_{1,2}$ and another corresponding to the remaining
eigenvalues. The first foliation has a cubic conjugate map and a quadratic
encoder, the second foliation has a linear conjugate map and a cubic
encoder. The results can be seen in figure \ref{fig:CarFollow-Auto}
for both the autonomous and the forced systems. The equation-driven
calculation breaks down before the amplitude reaches 0.15 - 0.18.
The data-driven calculation does not show this breakdown. 

\begin{figure}
\begin{centering}
\includegraphics[width=0.4\textwidth]{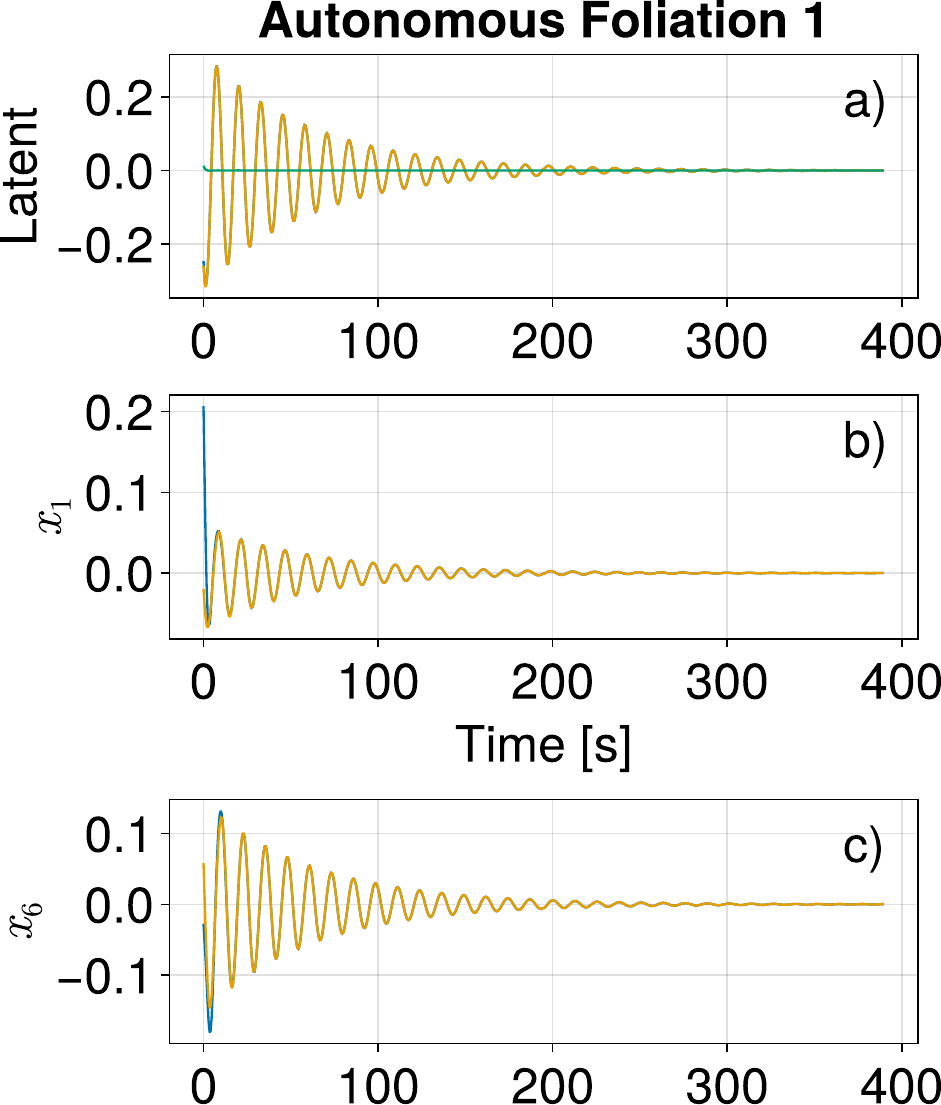}\includegraphics[width=0.4\textwidth]{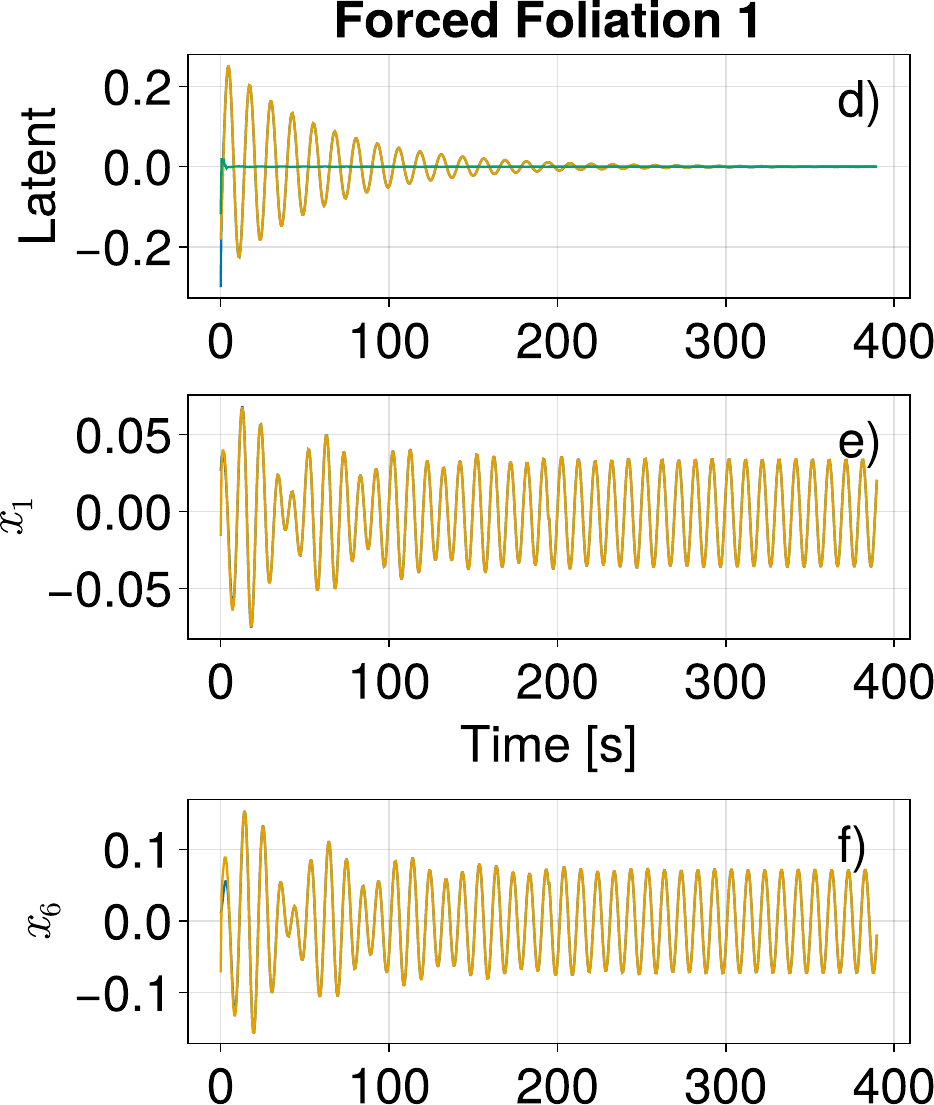}
\par\end{centering}
\caption{\label{fig:CarFollow-latent}Visualisation of the reduced order model
accuracy for model \eqref{eq:CarFollow-Model}. Interpretation of
the lines is the same as for figure \ref{fig:SP-latent}. The visualised
physical variables are the speed of the first car about the equilibrium
$x_{1}=v_{1}-v_{1}^{\star}$, and the headway of the first car about
the equilibrium $x_{6}=h_{1}-h_{1}^{\star}$. }

\end{figure}

Figure \ref{fig:CarFollow-latent} shows how the testing data fits
the model. The error is only noticeable for a fraction of a period
at the highest amplitudes. The full testing trajectory in physical
space is also reproduced with great accuracy, because the rest of
the dynamics decays very fast. Hence a two-dimensional reduced order
model provides a good representation of the whole system, which remains
the case with more cars.

\subsection{Clamped-clamped plate}

\begin{figure}
\begin{centering}
\includegraphics[width=0.99\textwidth]{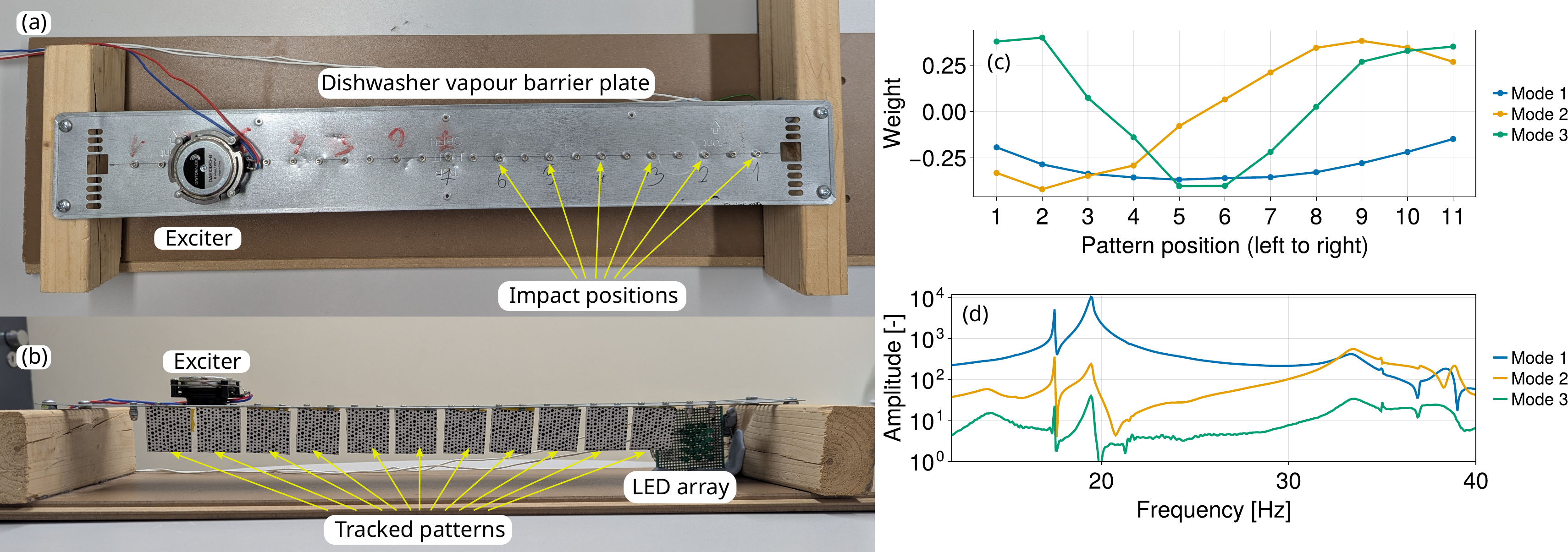}
\par\end{centering}
\caption{\label{fig:Plate-Rig}Experimental rig. a) View from the top. b) View
from the side. c) Mode shapes of the plate using PCA. d) Frequency
spectrum of the three vibration modes with forcing turned on.}
\end{figure}

This example shows the use of video data from a forced physical system.
The tested structure can be seen in figure \ref{fig:Plate-Rig}(a,b),
which consists of a mild steel plate mounted to two wooden blocks.
The forcing of the rig is produced by a Dayton Audio DAEX32Q-8 dynamic
exciter, which was driven at about 17.5 Hz. The same forcing signal
was also sent to an array of 20 LEDs, so that the phase of the forcing
can be recovered from the video. The plate was prepared with a series
of M3 screws. In order to better track motion a pattern of various
size dots were printed and attached to pairs of M3 screws using double
sided and electric tape. Two synchronised (Allied Vision, Alvium 1800
U-052M) cameras recorded the motion at 500 frames per second. The
spatial resolution of both cameras was $0.304$ mm/pixels. The rig
was then subjected to impact tests. The impacts were carried out by
a nylon mallet that hit the top of the M3 screws in 7 locations as
indicated in figure \ref{fig:Plate-Rig}(a). Two impact responses
were collected for each impact location, which resulted in 14 recorded
trajectories for the unforced and the forced system.

The motion from the video was extracted using digital image correlation
\cite{DIC2009}. Each pattern was treated as a rigid body and their
position at their centre and rotation was tracked. However, only the
vertical component of the motion was used, leading to 11 signals from
pairs of synchronised videos. Using principal component analysis,
three vibration modes were extracted from the 11 vertical displacements.
Any further vibration modes where not usable due to their noise content.
The three modes can be seen in figure \ref{fig:Plate-Rig}(c), which
strongly resemble beam vibration modes. The frequency spectrum of
an impact test on the forced system can be seen in figure \ref{fig:Plate-Rig}(d)
for each of the modes. Forcing at 17.5 Hz, the first mode at about
19.5 Hz and the second mode 33.2 Hz are the most noticeable. The rest
of the modes are mostly masked by noise, therefore we only focus on
the first two natural frequencies.

In order to reconstruct the state space, we used delay embedding of
length 29, which is roughly the number of samples within a forcing
period. This produced highly correlated data and therefore dynamic
mode decomposition was used to identify the first six dynamic modes
with the highest Euclidean norms, which also align with the slowest
decay rates. Three two-dimensional invariant foliations were fitted
to the data. The first two foliations had order-5 encoders of type
\eqref{eq:ENC-reducible} and order-5 conjugate maps. The third foliation
had a linear conjugate map and an order-5 locally defined encoder
of type \eqref{eq:ENC-local}.

\begin{figure}
\begin{centering}
\includegraphics[width=0.99\textwidth]{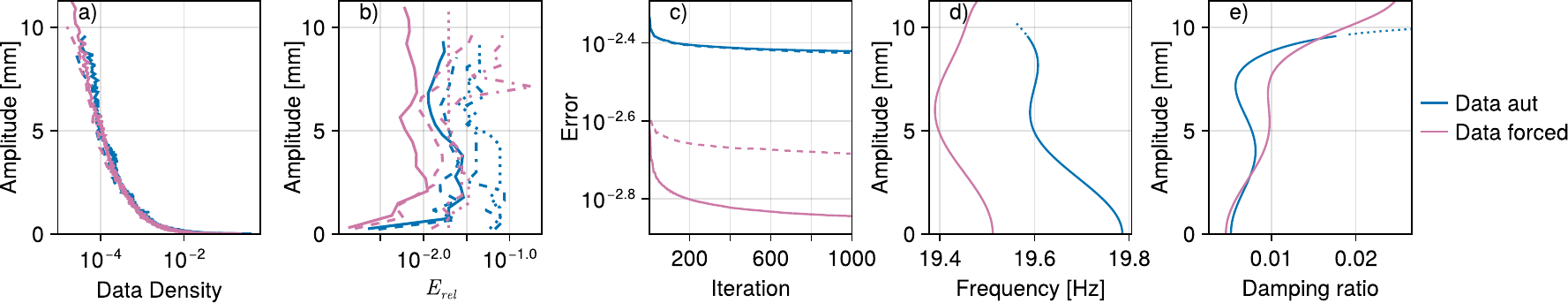}
\par\end{centering}
\caption{\label{fig:Plate-Res-1}First natural frequency of the plate in figure
\ref{fig:Plate-Rig}. The blue and purple lines represent the autonomous
and the forced system, respectively a) Data density: continuous lines
denote the training data dashed lines denote testing data. b) Relative
error \eqref{eq:Relative_Error}: continuous lines are mean relative
error of training data, dotted lines are maximum relative error of
training data, dashed lines are mean relative error of testing data,
dash-dotted lines are maximum relative error of testing data. c) Continuous
lines are the training errors, dashed lines are the testing errors
as functions of optimisation steps. d) Frequency backbone curve as
given by equation \eqref{eq:INST-freq}, dots mean extrapolation,
without training data. e) Damping ratio curve as given by equation
\eqref{eq:INST-damp}, dots mean extrapolation.}
\end{figure}
\begin{figure}
\begin{centering}
\includegraphics[width=0.99\textwidth]{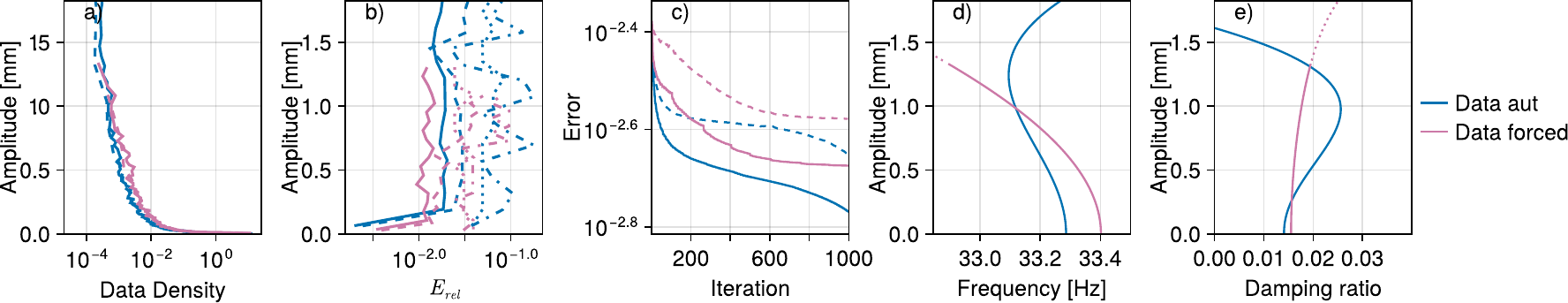}
\par\end{centering}
\caption{\label{fig:Plate-Res-2}Second natural frequency of the plate in figure
\ref{fig:Plate-Rig}. The legend is the same as in figure \eqref{fig:Plate-Res-1}.}

\end{figure}

The result for the first reduced order model can be seen in figure
\ref{fig:Plate-Res-1}, which is fully specified by the two backbone
curves \eqref{eq:INST-damp} and \eqref{eq:INST-freq} in figure \ref{fig:Plate-Res-1}(e)
and \ref{fig:Plate-Res-1}(d), respectively. The backbone curves for
the autonomous and forced case are similar, except that their frequency
is slightly shifted. The second vibration mode is resolved in figure
\ref{fig:Plate-Res-2}, with somewhat higher errors.

\begin{figure}
\begin{centering}
\includegraphics[width=0.49\textwidth]{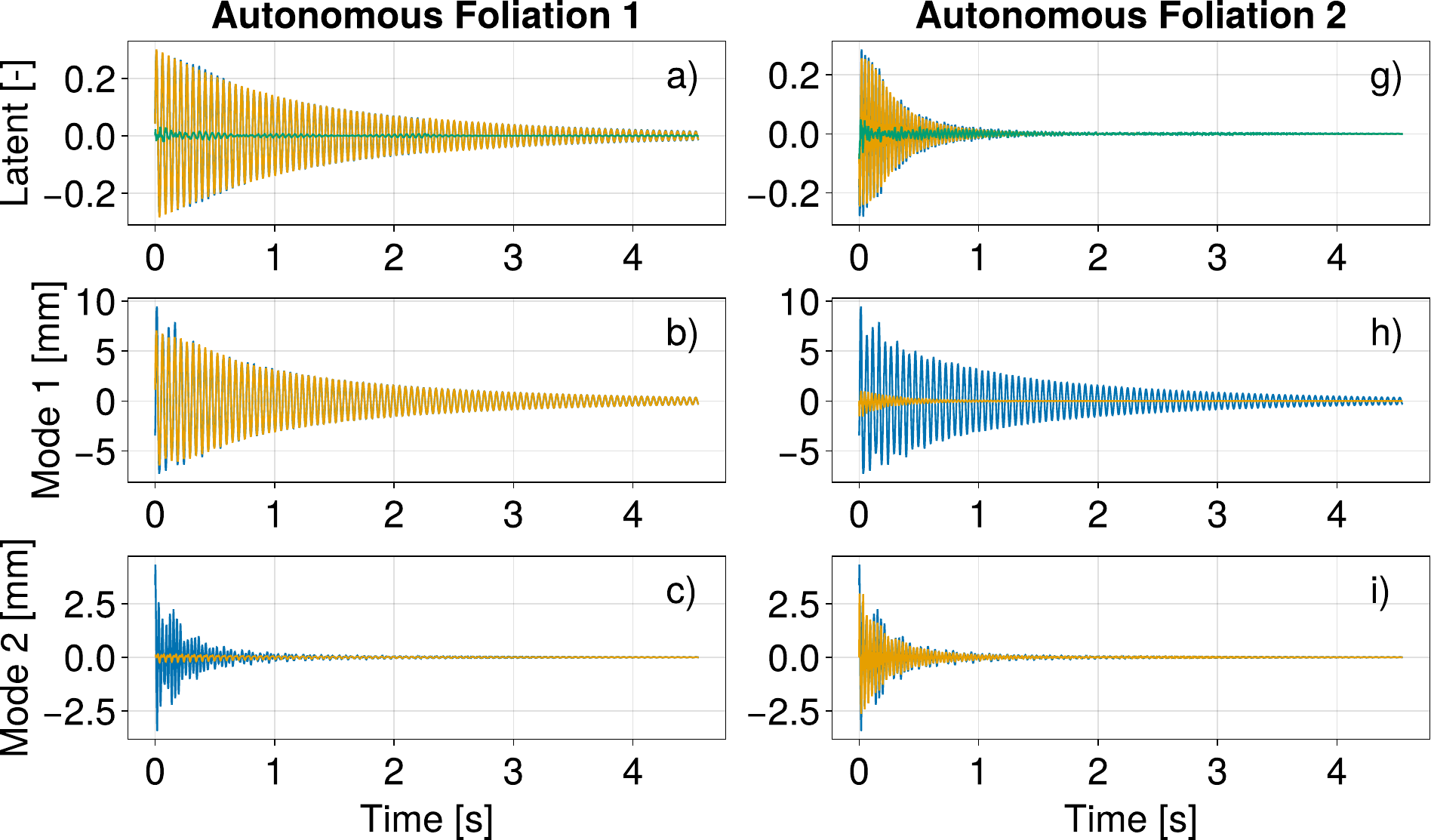}\includegraphics[width=0.49\textwidth]{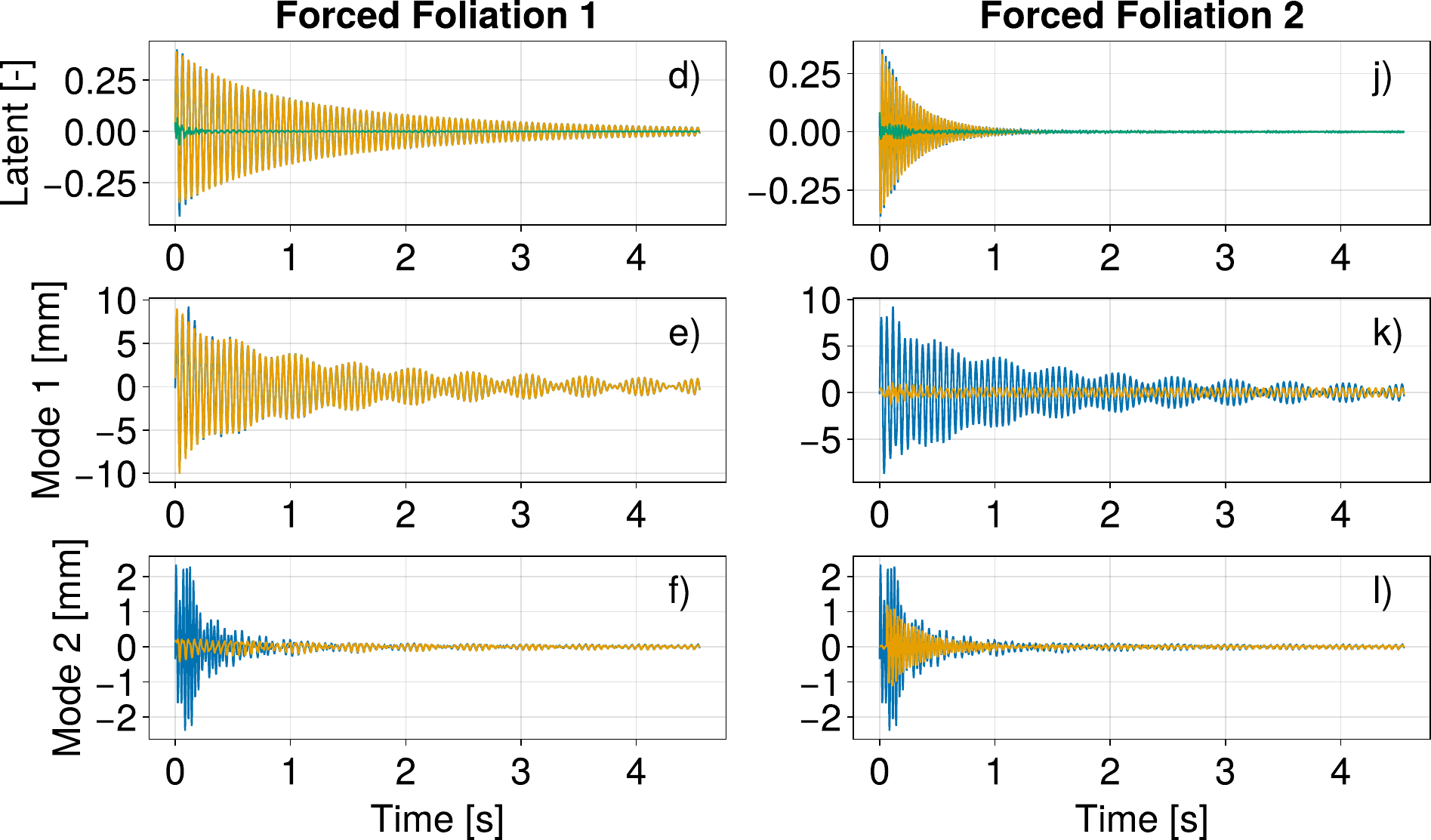}
\par\end{centering}
\caption{\label{fig:Plate-Latent}Visualisation of reduced order model accuracy
for the rig in figure \ref{fig:Plate-Rig}. Interpretation of the
lines is the same as for figure \ref{fig:SP-latent}.}
\end{figure}

The accuracy of the reduced order model can be also illustrated by
comparing the model prediction to a testing trajectory. In figures
\ref{fig:Plate-Latent}(a,d,g,j) we have plotted the encoded testing
data, the model prediction and their difference (green lines) in the
latent space. For the first natural frequency the accuracy is very
good, for the second natural frequency and for the autonomous system
the errors a more pronounced. However the vibration has the same phase
and amplitude as the encoded data, hence one can conclude that the
error is due to noise in the data that is not removed by the encoder.
It is also possible to map the model prediction back into the physical
space. In our case we are predicting the first two mode shapes (as
in figure \ref{fig:Plate-Rig}(c)). It can be seen that the model
of the first natural frequency almost fully reproduces the signal
in the first mode shape, apart from the initial transients represented
by higher frequency vibrations in figures \ref{fig:Plate-Latent}(b,e).
The model of the second natural frequency is less capable reproducing
the second mode (in figures \ref{fig:Plate-Latent}(i,l)) because
of noise in the data and stronger coupling to the first mode shape.

\section{Conclusions}

The paper has demonstrated that the invariant foliation architecture
is a practical and generally applicable tool to identify reduced order
models from data. The invariance and uniqueness of this architecture
makes sure that the produced model is both meaningful and reproducible.
The method is able to pick out multiple invariant subsets of the dynamics,
decomposing complex systems into smaller independent nonlinear components.
Using a single encoder reduces the number of parameters compared to
autoencoders that also need a decoder. We have found that in some
cases the identified models are able to extrapolate beyond the available
data. The accuracy of the identified model can be excellent, but also
depends on noise in the data.

\paragraph{Data availability and Software}

The data and computer code that reproduces the figures and calculations
can be found at \url{https://github.com/rsnumerics/InvariantModels.jl}
\cite{InvarianModels2025}.

\paragraph{Conflict of Interest}

The author declares that he has no conflict of interest. 

\paragraph{Funding}

No funding was received.

\printbibliography

\end{document}